%% file: main.tex
\documentclass[a4paper,fleqn]{cas-sc}
\usepackage{multirow}

\usepackage[numbers]{natbib}
\usepackage{lineno}
\usepackage{tikz}
\DeclareMathOperator{\sinc}{sinc}

\def\tsc#1{\csdef{#1}{\textsc{\lowercase{#1}}\xspace}}
\tsc{WGM}
\tsc{QE}
\tsc{EP}
\tsc{PMS}
\tsc{BEC}
\tsc{DE}

\begin{document}
\let\WriteBookmarks\relax
\def\floatpagepagefraction{1}
\def\textpagefraction{.001}

\shorttitle{Design principles for stable and generalizable data-driven discretizations}

\shortauthors{A.-A. Nasser and A. Adcroft}

\title [mode = title]{Design principles for stable and generalizable data-driven discretizations for solving linear hyperbolic conservation laws}

%

\author[1]{A.-A. Nasser}[type=editor,
                        auid=000,bioid=1,
                        orcid=0000-0001-6757-2650] 
\cormark[1] 
\fnmark[1] 
\ead{antoine.nasser@princeton.edu} 
\credit{Conceptualization of this study, Methodology, Software}

\author[1]{A. Adcroft}[type=editor,
                        auid=000,bioid=1,
                        orcid=0000-0001-9413-1017] 
\credit{Conceptualization of this study, Methodology, Software}
                        
\affiliation[1]{organization={Princeton University Atmospheric and Oceanic Sciences Program},
    city={Princeton},
    country={USA}}

\cortext[cor1]{Corresponding author}

\begin{abstract}
We investigate data-driven finite-volume discretizations of the linear advection equation in one dimension. Neural networks for use as numerical advection schemes are constructed adhering to first principles of numerical analysis, allowing us to examine how normalization, training data, and architectural choices influence stability, accuracy, and shape preservation. (i) We show that reconstruction based solely on cell averages leads to a multi-valued learning problem, explaining limited generalization when training data includes widely different curvature regimes. (ii) Numerical stability and good generalization can be achieved by enforcing semilinearity (Lin and Rood 1998) through local stencil-scale normalization, which ensures invariance under affine transformations of the inputs. (iii) A new data-driven flux limiter is introduced that outperforms the classical 'OSTVD3' (Arora and Roe, 1997) scheme in shape preservation by introducing mild antidiffusion in near-linear regimes, while higher-order reconstruction in non-monotonic regions provides limited benefit. (iv) We show that training on polynomial profiles yields stable, high-order accurate discretizations, with the polynomial degree controlling the formal order of accuracy. Together, these results illustrate how the representational, architectural, and training choices govern the stability and generalization of data-driven finite-volume schemes for linear advection.
\end{abstract}


\begin{highlights}
\item Exact advection from cell-averaged inputs is shape-specific, explaining poor generalization in ML-FV schemes.
\item Enforcing stencil-scale normalization yields semilinear schemes and prevents out-of-distribution behavior.
\item A data-driven flux limiter that outperforms OSTVD3 in shape preservation across smooth and discontinuous initial conditions.

\end{highlights}


\maketitle
\section{Introduction}

Accurate numerical treatment of transport processes is central to predictive simulations of geophysical flows \citep{balaji_are_2022}, astrophysical flows \citep{stone_athena_2008}, and engineering fluid dynamics \citep{ferziger_computational_2002}. Yet available computational resources limit the smallest physical scales that can be explicitly resolved. This motivates the investigation of machine-learning (ML) approaches to extend the effective resolution of numerical simulations. Despite recent progress, several outstanding challenges still prevent ML-based solvers from being routinely used in numerical models \citep{zanna_framework_2025}.

A central difficulty is the absence of a general framework for designing data-driven discretizations that both generalize and remain numerically stable. Previous studies \citep{bar-sinai_learning_2019, zhuang_learned_2021} have shown that enforcing mass conservation and at least first-order accuracy constitutes a minimal requirement to prevent rapid numerical blow-up. However, conservation alone is generally insufficient to avoid out-of-distribution behavior or instabilities over long time integrations \citep{bar-sinai_learning_2019, zhuang_learned_2021, morand_deep_2024, stevens_enhancement_2020}. For example, \citep{timofeyev_application_2025} embedded a neural network (NN)–based closure within a conservative flux-correction framework, achieving stable behavior and accurate predictions in scenarios not seen during training. Part of this performance, however, depends on the flux limiting strategy, which may be applied a posteriori or incorporated during training. Alternatively, \citep{bar-sinai_learning_2019} account for NN errors during training by introducing a multi-step loss that reduces the accumulation of space–time discretization errors and improves stability. Such approaches, however, do not provide a priori guarantees, and the number of steps required to ensure stability remains unclear \citep{kochkov_machine_2021}. \citep{alieva_toward_2023} argue that selecting turbulent training data from statistically steady states enables stable long-time integration, as point-wise accuracy can be traded for long-time behavior under ergodic arguments. However, their results apply to global operators rather than local schemes, and it remains unclear whether this equivalence holds at the scale of local stencils. In contrast, classical numerical analysis interprets spatial and temporal errors through truncation terms in the discretization \citep{lax_systems_1960}, providing explicit stability conditions. This contrast raises the question of whether stability can instead be enforced within a unified numerical framework. In addition, \citep{zhuang_learned_2021} showed that semilinear data-driven schemes \citep{lin_multidimensional_1996}, which preserve linear tracer correlations, significantly enhance numerical stability. From a machine-learning perspective, this property corresponds to rescaling both inputs and outputs, thereby mitigating distribution shifts in the training data \citep{beucler_climate-invariant_2024, kim_reversible_2021}. Nevertheless, the data-driven discretization proposed by \citep{zhuang_learned_2021} remains prone to overfitting, with solutions tending toward shapes observed during training. Hence, several barriers remain in achieving robust generalization for learned discretizations.

These challenges are compounded by the fact that the design of traditional advection schemes remains an open problem. A well-known difficulty is that formally high-order discretizations often produce spurious oscillations (or wiggles) near shocks or steep gradients, which must be controlled through flux or slope limiters \citep{leveque_numerical_1992}. As originally formulated \citep{colella_piecewise_1984,arora_well-behaved_1997,harten_high_1983}, limiters have the property that the truncation error is first-order accurate at all extrema, regardless of the accuracy of the underlying high-order method. This limitation is solved, for instance, in \citep{colella_limiter_2008} where a high-order reconstruction is introduced near smooth extrema which significantly improves the accuracy of the solution. In practice, numerous limiter formulations have been proposed, yet no single method performs optimally across all problems \citep{zhang_review_2015}. Machine-learning approaches may offer new insights into the design of improved flux limiters. For instance, \citep{nguyen-fotiadis_machine_2022, nguyen-fotiadis_probabilistic_2025} argue that spring-like formulations of limiters in the Sweby diagram improve shape preservation for Burgers’ equation. In contrast, \citep{huang_learning_2025} showed that a simple combination of the Superbee and Minmod limiters \citep{roe_asymptotic_1983} can surpass classical limiters across a range of problems with shocks and discontinuities.

To address these questions, we investigate strategies for improving the numerical stability and generalization of data-driven discretizations for transport problems in a one-dimensional periodic domain. First, we cover the numerical framework required for constructing stable data-driven discretizations (Sections 1.1-1.3). Next, we demonstrate that reconstructing solutions from cell-averaged inputs leads to a multi-valued learning problem, which limits the generalization ability of data-driven discretizations (Section 2). We then introduce a data-driven three-point flux limiter that improves shape preservation during simulated transport and outperforms classical TVD-based formulations \citep{harten_high_1983} (Section 3). Finally, Sections 4–5 summarize the main findings and discuss their implications.

\subsection{Linear one-dimensional advection problem}

In this study, we focus on the solution of the linear scalar advection equation
\begin{eqnarray} \label{eq:adv}
d_t q = \partial_t q + u \partial_x q = 0, 
\end{eqnarray}
where $u$ is a constant velocity and for which the general solution is $q(x,t)=q(x-ut,0)$. In many applications, the conservation laws are written in flux-form, i.e.
\begin{eqnarray}
\partial_t q + \partial_x (uq) = 0.
\end{eqnarray}
For simplicity, we assume $u>0$; the case $u<0$ follows by symmetry with respect to the cell interface. Let $\tau$ and $\Delta$ denote the time step and the cell width, respectively. The discrete solution at time $t^n = n\tau$ and spatial location $x_i = i\Delta$ is denoted by $q_i^n$.

\subsection{Finite-volume approximation}
Following \citep{eymard_finite_2000}, the finite-volume (FV) formulation represents the discrete solution with volume average values given by
\begin{eqnarray} \label{eq:cavg}
q_i^n = \frac{1}{\Delta} \int_{x_{i-1/2}}^{x_{i+1/2}} q(x,t^n) dx,
\end{eqnarray}
where $x_{i\pm 1/2}$ denote the cell interfaces of cell $i$. The one-step approach \citep{lax_systems_1960} leads to the update
\begin{eqnarray} \label{eq:discrete}
q_i^{n+1} = q_i^{n} -\frac{\tau}{\Delta} \left( F_{i+1/2}^n - F_{i-1/2}^n \right),
\end{eqnarray}
where $F_{i+1/2}^n$ is the numerical flux through the right cell interface. The exact time-averaged flux can be written as
\begin{eqnarray}
F_{i+1/2}^n = \frac{1}{\tau} \int_{t^n}^{t^{n+1}} u q(x_{i+1/2},t) dt,
\end{eqnarray}
or, using the characteristic solution above, equivalently,
\begin{eqnarray} \label{eq:cflux}
F_{i+1/2}^n = \frac{u}{\nu\Delta} \int_{x_{i+1/2}-\nu \Delta}^{x_{i+1/2}} q(x,t^n) dx
\end{eqnarray}
where $\nu = u\tau / \Delta$ is the Courant number. This expression corresponds to the spatial average of $q$ over an interval of length $\nu\Delta$ upstream of the interface. 

\subsection{Coupled one-step approach}\label{sect:one-step}
Following \cite{daru_high_2004}, discrete advection schemes are constructed by successively correcting the error terms in the modified equation arising from the one-step approach \citep{lax_systems_1960}, yielding schemes that are accurate in both space and time (referred to as 'OS' schemes). In this framework, the discrete flux \eqref{eq:cflux} can be written in the form of a Lax–Wendroff scheme modified by an accuracy function $\Phi$,
\begin{equation} \label{eq:cphi}
     F^{n}_{i+1/2} =  q^n_i + \frac{1-\nu}{2} \Phi_{i+1/2} \left( q^n_{i+1} -  q^n_{i} \right).
\end{equation}
Setting $\Phi = 1$ recovers the second-order Lax–Wendroff scheme ('OS2'), while $\Phi = 0$ yields the first-order upwind scheme ('OS1'). A third-order accurate scheme ('OS3') is obtained by defining
\begin{eqnarray}
    \Phi^{\text{3}}_{i+1/2} = 1 - \frac{1+\nu}{3}(1-r_{i+1/2})
\end{eqnarray}
where the slope ratio $r_{i+1/2} = (q_i - q_{i-1})/(q_{i+1} - q_i)$ is a measure of local curvature (of the discrete values, $q_i$): $r=1$ for a straight line, $r\to 0$ or $r\to\infty$ for monotonic sharp gradients, and $r<0$ in non-monotonic regions.  By successively correcting truncation errors, arbitrarily high-order schemes can be constructed by defining $\Phi$, with explicit expressions provided in \cite{daru_high_2004, del_pino_arbitrary_2006}. 

\section{Numerics and generalization of data-driven discretizations}
\subsection{Physics-informed neural networks for solving advection}
\subsubsection{Hybrid model}
\begin{figure}
    \centering
    \includegraphics[width=\linewidth]{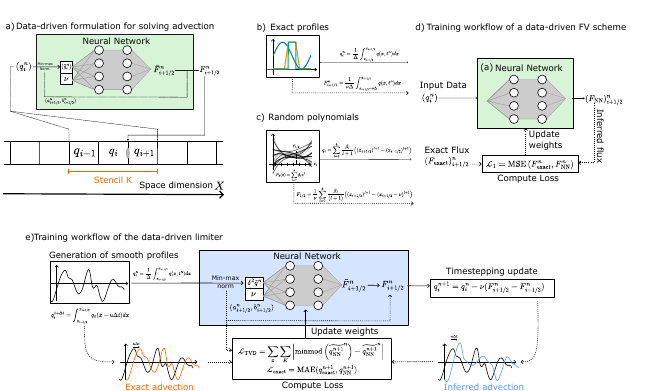}
    \caption{Illustration of the problem setup and architectures of the data-driven discretizations. (a) Data-driven formulation for the advection problem. (b–c) Generation of training data using analytical profiles (b) and random polynomials (c). (d) Under local input–output normalization (as in panel a), the network is trained on exact pairs of FV cell-averages and numerical fluxes. (e) Training workflow of the data-driven limiter: random harmonic profiles are generated, and the loss with an added TVD penalty term is computed between the exact FV solution and the neural network predictions.}
    \label{fig:method}
\end{figure}

On the basis that the OS schemes are accurate and stable, we now investigate whether neural networks can be trained to work in the same context, inheriting the stability, but with better accuracy. Here, the neural network predicts the numerical fluxes in \eqref{eq:discrete} at time $t^n$ to advance the solution to $t^{n+1}$, using as inputs the local grid values $(q_{i+k}^n)_{k\in K}$ (Fig.~\ref{fig:method}a). The network weights are learned by minimizing the mean square errors (MSE), averaged over the domain, between the prediction $F_{\mathrm{NN}}^n$ and the exact finite-volume flux given from Eq.\ref{eq:cflux} (see Fig.~\ref{fig:method}d), giving
\begin{eqnarray}
\mathcal{L}_1 = \sum_x \left(F_{\mathrm{exact}}^n - F_{\mathrm{NN}}^n\right)^2.
\end{eqnarray}

Global conservation of mass is a minimal requirement for stable conventional advection schemes and is ensured if the numerical flux at a cell interface is uniquely defined, independent of whether it is evaluated from the left or from the right cell. Enforcing this property in a machine learning based approach is therefore essential and is ensured by predicting a unique flux. 

\subsubsection{Local semilinearity invariant} \label{sect:minmax}
The transport equation Eq.~\ref{eq:adv} is semilinear invariant, meaning that affine relationships between tracers are preserved during advection, a property that discrete schemes should also satisfy \citep{lin_multidimensional_1996}. \citep{zhuang_learned_2021} argued that satisfying semilinearity improved numerical stability which they achieved by applying global normalization of the solution prior to neural-network inference, followed by denormalization of the predicted flux. From a machine-learning perspective, this would require sufficiently dense sampling of the hypercube $[0,1]^N$, with $N$ the number of grid points, to avoid out-of-distribution inputs, a requirement that becomes prohibitive as $N$ increases and is therefore intractable in practice. Here, we suggest instead to apply normalization locally at the stencil scale, thereby ensuring that semilinearity holds at the level of the numerical flux. Local normalization provides three main advantages: it preserves tracer correlations at the scale of a few grid cells; it guarantees preservation of constant states, a consistency condition required for stability \citep{bar-sinai_learning_2019}; and it defines a stronger invariance property than global normalization, which ensures no out-of-distribution behavior and simplifies the training dataset \citep{kim_reversible_2021}.

At the discrete level, the flux is equivariant under affine transformations of the inputs, i.e.
\begin{eqnarray}\label{eq:renorm}
F_{i+1/2}^n = a_{i+1/2}^n \tilde F_{i+1/2}^n + b_{i+1/2}^n
\end{eqnarray}
where the discrete normalization factors $(a_{i+1/2}^n,b_{i+1/2}^n)$ are defined over the input stencil $K$ used to infer $F_{i+1/2}^n$, tilde refers to the normalized solution related by $q = a \tilde q + b$ and $\tilde F$ is the flux associated with $\tilde q$, written
\begin{eqnarray}\label{eq:cfluxnorm}
\tilde F_{i+1/2}^n = \frac{u}{\nu\Delta} \int_{x_{i+1/2}-\nu\Delta}^{x_{i+1/2}} \tilde q(x,t^n)dx.
\end{eqnarray}

Several scaling approaches are possible to enforce local semilinear property in the neural network, here, we adopt the min–max normalization rule:
\begin{eqnarray}
a_{i+1/2}^n &=& \max_{k\in K}(q_{i+k}^n) - \min_{k\in K}(q_{i+k}^n) \\
b_{i+1/2}^n &=& \min_{k\in K}(q_{i+k}^n)
\end{eqnarray}
and define normalized inputs
\begin{eqnarray}
\tilde q_i^n = \frac{q_i^n - b_{i+1/2}^n}{a_{i+1/2}^n}.
\end{eqnarray}
Hence, a neural network takes normalized data $(q_i^n)_i$ as input to predict the normalized flux $\tilde F_{i+1/2}^n$, which gives the physical flux $F_{i+1/2}^n$ using \eqref{eq:renorm}. To prevent singularities when $a_{i+1/2}^n \to 0$, the scaling factor is clipped to machine precision.

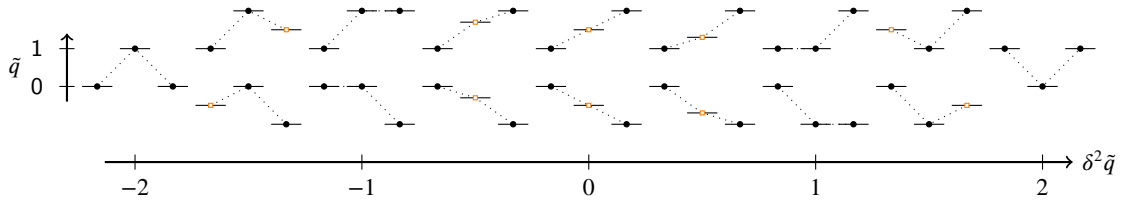
\begin{figure}
    \centering   
    \input{figures/fig_stencil_curv2}
    \caption{Possible configurations of a three-point stencil in normalized phase space, sorted according to the curvature of the gridded data. The horizontal axis represents the local curvature of the normalized field, $\tilde q_{i-1}-2\tilde q_i+\tilde q_{i+1}$. Black dots denote values equal to 0 or 1 depending on their vertical position, whereas intermediate values are shown in orange. Configurations located in the upper (lower) part of the figure correspond to positive (negative) slopes, defined by $\tilde q_{i+1}-\tilde q_{i-1}$. The symmetric-extrema case defines the bounds of the admissible configurations.}
    \label{fig:K3}
\end{figure}

Now consider the level of representation of a discrete scheme using normalized inputs in function of its stencil size. For $K=1$, both the normalized value, $\tilde q_i^n$, and the normalized flux, $\tilde F_{i+1/2}^n$, vanish, making the normalized scheme degenerate. For $K=2$, three configurations exist depending on $(q_i^n,q_{i+1}^n)$: a plateau $q_i^n=q_{i+1}^n$, a positive jump $q_i^n<q_{i+1}^n$, and a negative jump $q_i^n>q_{i+1}^n$. After normalization, these reduce to $(0,0)$ and $(0,1)$, with the sign of the jump encoded in $a$. Therefore, the smallest stencil that yields expressive behavior is $K=3$, since one free continuous parameter remains after normalization. In this case, six distinct permutations arise, characterized by the mean slope, $\tilde q_{i+1}-\tilde q_{i-1}$, and the position of the intermediate value in the stencil, allowing the discrete representation of arbitrary functions up to quadratic order (see Fig.~\ref{fig:K3}). This interpretation extends to larger stencils, by introducing additional free parameters which enables the representation of higher-order curvatures. We note that by virtue of local normalization, if the training data includes the samples at the transitions between permutations then we can be sure the input features at inference will always be within the range of the training data. We now investigate data-driven discretizations based on a three-point stencil ($K=3$). Throughout this study, the local curvature of the normalized field is represented by the discrete second derivative, $\delta^2 \tilde q = \tilde q_{i-1}-2\tilde q_i+\tilde q_{i+1}$

\subsection{Training with smooth solutions recovers traditional accurate schemes} \label{sect:multivalued}

Following \eqref{eq:cavg} and \eqref{eq:cflux}, a given continuous solution uniquely determines both the cell-averaged values and the associated numerical flux. However, the converse is not true since a discrete representation may correspond to multiple underlying profiles, raising the question of whether a discretization using cell-averaged inputs can generalize well across distinct solutions.

\begin{figure}
    \centering
    \hspace*{-3cm}
    \includegraphics[width=\linewidth, trim={0 0 18cm 0}]{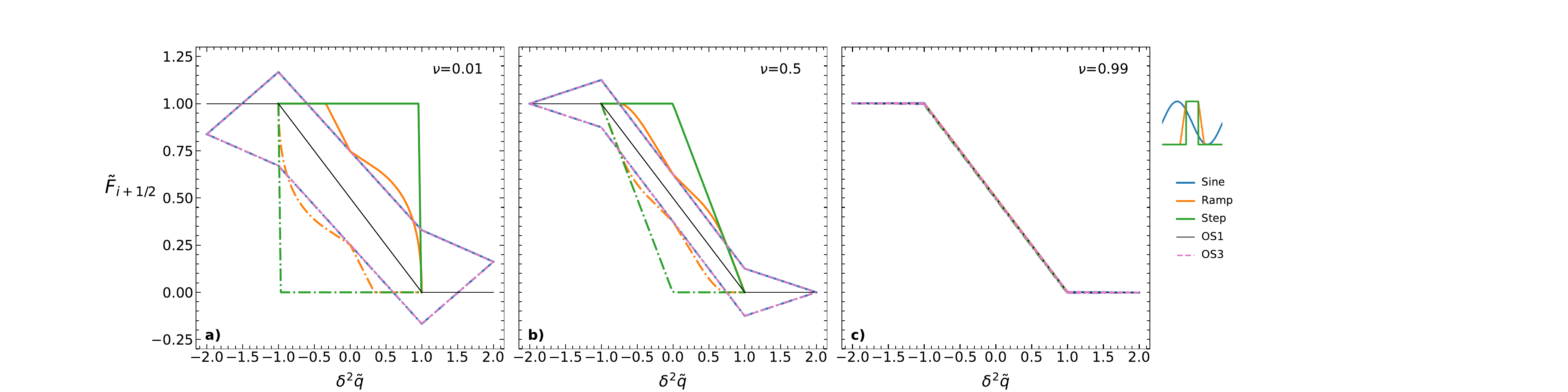}
    \caption{Flux-curvature diagrams for perfect advection of three initial conditions for a three-point stencil. Each color denotes a shape specified in the legend. The normalized numerical flux $\tilde F_{i+1/2}$ is plotted in the $y$-axis against the local curvature of the normalized gridded field, $\tilde q_{i-1}-2\tilde q_{i}+\tilde q_{i+1}$ along $x$-axis. The curvature range corresponds to the same stencil values shown in Fig.~\ref{fig:K3}, where the flux for positive slopes (top part of Fig.~\ref{fig:K3}) are shown as solid lines, and those for negative slopes (bottom) as dotted dashed lines. Each panel shows the flux at a CFL number of (a) 0.01, (b) 0.5, and (c) 0.99.}
    \label{fig:learnable}
\end{figure}

To address this question, we evaluate the exact finite-volume flux, within the phase space of a normalized three-point stencil, of three analytical profiles commonly used for training: an abrupt step, a linear ramp, and a sine wave (Fig.~\ref{fig:learnable}). The cell-averaged values and the flux are calculated from exact solution and sampled by incrementally shifting the profile across the full stencil length (the sampling procedure is detailed in Appendix~\ref{sect:appendix-A}). Each profile induces a trajectory in phase space which, if learned, allows simulating perfect advection of the profile. From left to right, panels a, b, and c correspond to CFL numbers of 0.01, 0.5, and 0.99, respectively.

Figs.~\ref{fig:learnable}a,b show that the three profiles produce distinct mappings in phase space. This behavior is consistent with Eq.~\eqref{eq:cflux}, where the numerical flux approaches the mid-point value as the CFL number tends toward zero. In contrast, Fig.~\ref{fig:learnable}c shows that in the limit of pure advection, when the CFL number approaches one, all mappings collapse onto the first-order upwind scheme (black line in Fig.~\ref{fig:learnable}). A consistent feature across Figs.~\ref{fig:learnable}a–c is that the sine-wave trajectory (blue line) closely aligns with the third-order reconstruction OS3 (dashed pink line). This agreement reflects that, for smooth profiles (i.e. the profile curvature is well resolved at a given Courant number), the third-order scheme recovers the Taylor expansion of the sine function up to $O(\Delta^3)$. Sensitivity tests using other smooth test functions produce similar trajectories (not shown), indicating that the numerical flux associated with smooth solutions is already accurately reproduced by classical high-order schemes.

Fig.~\ref{fig:learnable} demonstrates that discrete representations using cell-averages leads to a multi-valued problem when training solutions range from discontinuous to smooth shapes. This multi-valued mapping reflects the fact that cell averages alone do not uniquely determine the sub-grid profile, so no single deterministic mapping from a fixed stencil of cell averages to a numerical flux can exactly reproduce the advection of profiles with differing regularity.

\subsection{Impact of training data on stability and generalization}\label{sect:2.3}
\begin{figure}
    \centering
    \includegraphics[width=0.6\linewidth]{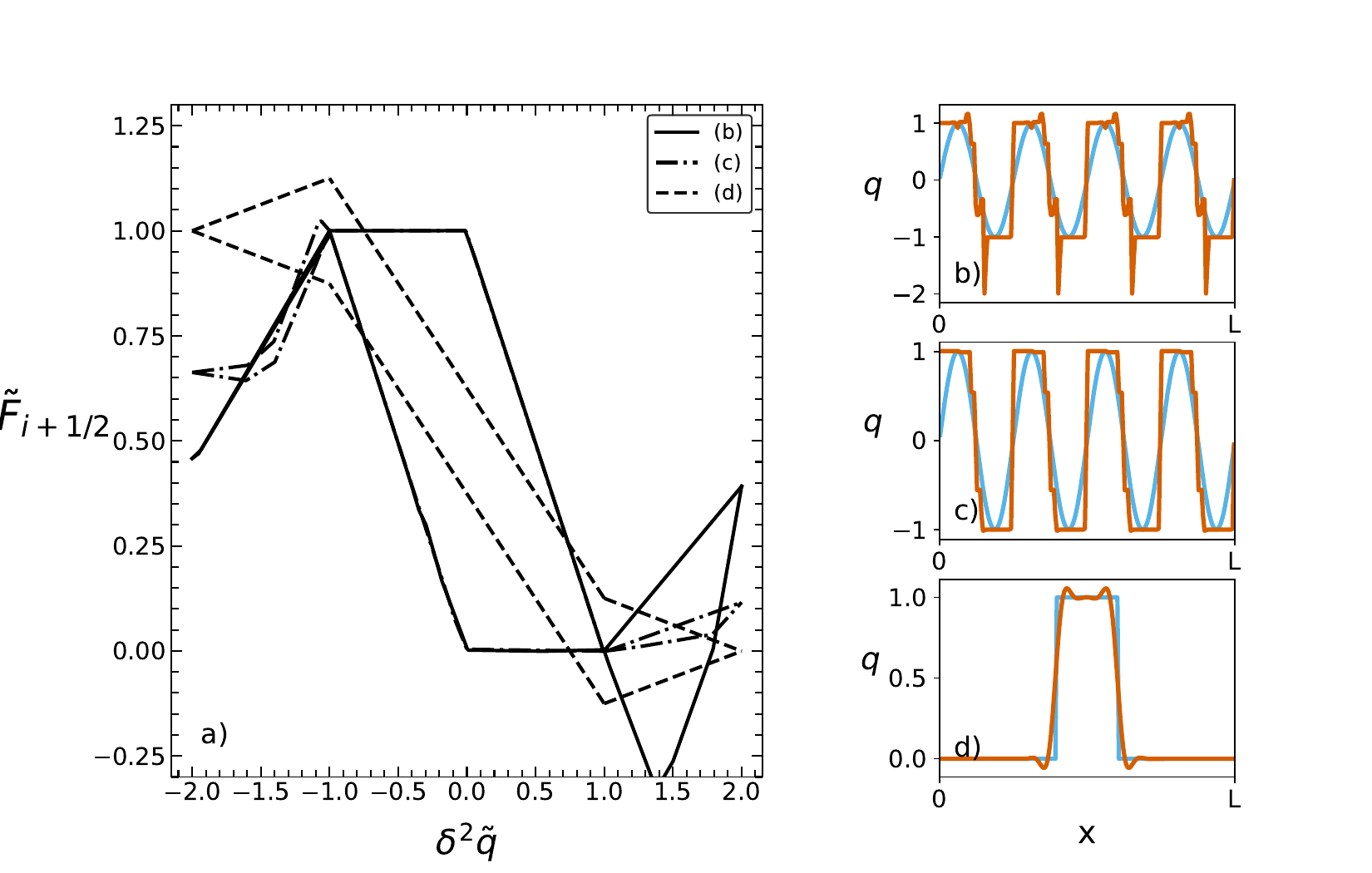}
    \caption{Flux-curvature diagram for data-driven discretizations trained on specific solutions. Panel (a): phase space of neural networks trained to predict the flux for a step profile (two random seeds (b,c)) and for a sine profile (d); axes are identical to those in Fig. \ref{fig:learnable}. Numerical simulations with $\nu=0.5$ show numerical instability after a few time steps (b) and distorted solutions after three full domain cycles (c,d). The blue line denotes the initial condition and the orange line the model response; the y-axis shows solution $q$ and the x-axis the spatial domain [0, L].}
    \label{fig:learned}
\end{figure}

A recurring challenge is that trained models may depend on random seeds used in training, which can then manifest as undetermined behavior during integration of a solution \citep{bar-sinai_learning_2019}. Figure~\ref{fig:learned} illustrates this issue by comparing the advection of initial conditions that lie outside the training distribution for different neural networks. The hyperparameters and training data are described in Appendices~\ref{sect:appendix-B} and \ref{sect:appendix-C}, respectively.

In Figs.~\ref{fig:learned}b–c, the advection of a sinusoidal initial condition is shown for two neural networks trained on exact step solutions and initialized with different random weights. The corresponding learned mappings in phase space are displayed in Fig.~\ref{fig:learned}a. In Fig.~\ref{fig:learned}b, extrapolation of the network in non-monotonic regions leads to numerical instability, causing the solution to blow up after only a few time steps. In contrast, in Fig.~\ref{fig:learned}c the integration remains stable but exhibits strong overfitting: the initially smooth profile evolves into a square-like wave after one domain cycle, consistent with Fig.~4 in \citep{zhuang_learned_2021}. The numerical instability observed in Fig.~\ref{fig:learned}b can be explained by the choice of training data, since step-like profiles leave non-monotonic regions of the mappable phase space unconstrained. By contrast, training on harmonic functions spans the full phase space, which ensures numerically stable advection even for out-of-sample initial conditions, as illustrated in Fig.~\ref{fig:learned}d for an initial square profile.

\subsection{Towards high-order accurate and stable data-driven discretizations}\label{sect:2.4}
As the highest-order curvature representable on a stencil of size $K$ is of degree $K-1$, we introduce a framework for constructing high-order one-step data-driven (OSDD) schemes. The central idea is to learn accurate reconstruction with neural networks using polynomials defined at the stencil scale and of maximum degree $K-1$ (as illustrated in Fig.~\ref{fig:method}c). In theory, training on local polynomials should yield a formal order of accuracy equal to the polynomial degree \citep{daru_high_2004}. The hyperparameters are gathered in Appendix~\ref{sect:appendix-B} and the procedure used to generate the random polynomial training data is described in Appendix~\ref{sect:appendix-C}. We show that this framework provides one-step high-order schemes as summarized in Table~\ref{tab:accpoly}, where the OSDD index equals the order of accuracy $\mathcal{O}(\Delta^K)$. The OSDD approach remains accurate, at the appropriate order, provided that neural-network approximation errors are adequately controlled (e.g., \citep{wang_multi-stage_2024}). We found that computational inaccuracies began to interfere with the empirical estimation of $L_2$ order for $K>5$, likely due to the reduced dynamic range in the $L_2$ error; the network software introduces errors of $O(10^{-12})$. While learning discretizations that can be derived by analysis seems redundant, applying the data-driven methodology to multidimensional and multi-step settings does have the advantage of alleviating the need for cumbersome analytical derivations \citep{del_pino_arbitrary_2006}.

\begin{table}[]
    \centering
    \caption{Transport of a sum of 4 sine waves: $L_2$ error and order of accuracy for the one-step data-driven (OSDD) schemes.}
    \begin{tabular}{llll}
         Method & Number of grid points & $L_2$ error & $L_2$ order \\
         OSDD3 &  160  & 2.738 $\times 10^{-3}$ &     \\
                &  320  & 3.469 $\times 10^{-4}$ & 2.98 \\
                &  640  & 4.352 $\times 10^{-5}$ & 3.00 \\
                & 1280  & 5.445 $\times 10^{-6}$ & 3.00 \\
                & 2560  & 6.807 $\times 10^{-7}$ & 3.00 \\
            &&&\\
            OSDD4 &  160  & 5.820 $\times 10^{-4}$ &     \\
                &  320  & 3.706 $\times 10^{-5}$ & 3.97 \\
                &  640  & 2.327 $\times 10^{-6}$ & 3.99 \\
                & 1280  & 1.456 $\times 10^{-7}$ & 4.00 \\
                & 2560  & 9.103 $\times 10^{-9}$ & 4.00 \\
            &&&\\
            OSDD5 &  160  & 8.187 $\times 10^{-5}$ &     \\
                &  320  & 2.617 $\times 10^{-6}$ & 4.97 \\
                &  640  & 8.223 $\times 10^{-8}$ & 4.99 \\
                & 1280  & 2.573 $\times 10^{-9}$ & 5.00 \\
                & 2560  & 8.044 $\times 10^{-11}$ & 5.00 \\
    \end{tabular}
    \label{tab:accpoly}
\end{table}

\section{Learning enhanced shape-preserving flux limiter with data-driven advection}\label{sec:ddl}
\subsection{Background}
Among the monotonicity-preserving schemes, the Total Variation Diminishing (TVD) conditions \citep{harten_high_1983} are generally considered well suited for the preservation of sharp features but also considered too diffusive near extrema due to the limitation to first-order upwinding \cite{daru_high_2004}. The TVD conditions ensures that an unlimited scheme is monotonicity-preserving, if (assuming $\Phi=0$ for $r<0$)
\begin{eqnarray} \label{eq:harten}
\begin{cases}
0 \leq \Phi_{i+1/2} \leq \frac{2r_{i+1/2}}{\nu},  \\
0 \leq \Phi_{i+1/2} \leq \frac{2}{1-\nu}.
\end{cases}  
\end{eqnarray}
Following \citep{arora_well-behaved_1997}, enforcing the TVD constraints on the third-order scheme $\Phi^{3}$ yields a method that is formally third-order accurate almost everywhere and avoids the generation of spurious oscillations near sharp features (denoted OSTVD3):
\begin{eqnarray}
    \Phi^{\text{OSTVD3}}_{i+1/2}  = \max{\left(\min{\left(\Phi^{\text{3}}_{i+1/2},\frac{2r_{i+1/2} }{\nu},\frac{2}{1-\nu}\right)},0\right)}.
\end{eqnarray}
If instead the accurate transport of fronts and sharp gradients is the primary objective, the second-order shock-capturing limiter (so-called “Superbee” limiter), introduced in \citep{roe_asymptotic_1983}, is often preferred. It writes
\begin{eqnarray}
    \Phi^{\text{Superbee}}_{i+1/2} = \max{\left(
    \min{\left(\frac{2r_{i+1/2} }{\nu},1 \right)},
    \min{\left(r_{i+1/2},\frac{2}{1-\nu}\right)},
    0\right)}.
\end{eqnarray}
Hence, the choice of limiter remains application-dependent, and designing a formulation that gives accurate numerical solutions while minimizing distortions continues to be an open question.

\subsection{Data-driven approach to flux-limiter design} \label{sect:ddfl}

It follows from the previous section that training a discretization solely for formal accuracy yields solutions specific to the training dataset and is therefore prone to overfitting. An alternative approach to design a generalizable discretization is to train the network using both accuracy and monotonicity criteria to promote shape-preservation. By doing this we explore whether a flux-limiter formulation, applicable to linear conservation laws, can be learned within machine-learning framework. We now present this construction in the case of a three-point stencil.

\subsubsection{Neural network architecture}\label{sect:mirror}
We consider two parameterizations of the local stencil (Fig.~\ref{fig:K3}) as inputs to the neural network, together with the CFL number. The slope ratio $r$ is a natural choice in flux-limiter formulations, however, it is poorly suited from machine-learning perspective, as it becomes singular at extrema where its denominator vanishes, making it more prone to out-of-distribution behavior. A desirable feature of $r$-based formulations is their invariance to the sign of the gradient, where positive and negative slopes are treated identically. To retain this symmetry while avoiding singular behavior, we instead use the locally normalized curvature of the discrete field, $\delta^2 \tilde q$ (Section~\ref{sect:minmax}), as input. Symmetry is enforced by introducing a signed curvature: for negative slopes, $\tilde q_{i+1} - \tilde q_{i-1} < 0$, we apply the transformation $\delta^2 \tilde q \mapsto -\delta^2 \tilde q$, ensuring an identical representation of up- and down-gradient configurations. This symmetry implies the corresponding mapping of the flux, $\tilde F \mapsto 1 - \tilde F$. The hyperparameters of the neural network are gathered in Appendix~\ref{sect:appendix-B}.

The training workflow of the data-driven limiter is summarized in Fig.~\ref{fig:method}e. The central idea is to combine an accuracy loss with a penalty term that promotes monotonicity preservation. Minimizing the global mean-absolute error (MAE) between the neural-network prediction, $q_{\mathrm{NN}}^{n+1}$, and the exact solution, $q_{\mathrm{exact}}^{n+1}$, drives the network toward the unlimited third-order scheme, $\Phi^3$. This loss is defined as
\begin{eqnarray}
   \mathcal{L}_\text{exact} = \sum_x \bigg|q_{\mathrm{exact}}^{n+1} - q_{\mathrm{NN}}^{n+1} \bigg|.
\end{eqnarray}
To remove the emergence of spurious oscillations near sharp gradients, we introduce a TVD penalty inspired by the Local Discrete Maximum principle \citep{lipnikov_minimal_2012}, which discourages reconstructed values from leaving the initial range of the reconstruction stencil. It writes
\begin{eqnarray}
   \mathcal{L}_\text{TVD} = \sum_x\sum_K \bigg|\mathrm{minmod}\left(\widetilde{q^{n+1}_{\mathrm{NN}}}^n \right) - \widetilde{q^{n+1}_{\mathrm{NN}}}^n \bigg|,
\end{eqnarray}
with $\widetilde{q^{n+1}_i}^n = (q_i^{n+1} - b_{i+1/2}^n)/{a_{i+1/2}^n}$, and $\mathrm{minmod}(x)=\max(\min(x,1),0)$. The two losses are combined using an empirically determined penalty parameter $\lambda=0.5$:
\begin{eqnarray}
\mathcal{L}_{2} = (1-\lambda) \mathcal{L}_\text{exact} + \lambda \mathcal{L}_\text{TVD}\;.
\end{eqnarray}

As local modifications induced by a flux limiter influence the global shape preservation of the solution, we train on harmonic functions defined over the full domain to include these non-local effects in the loss. Data generation is detailed in Appendix~\ref{sect:appendix-C}.

\subsubsection{Sensitivity to weight initialization}

Numerical experiments indicate that the final performance of the learned limiter depends on the random initialization of the network weights. To assess this sensitivity, we train a set of 20 models with different random seeds and rank them using a shape-preservation metric. This metric is defined as the mean $L^1$ error after one full cycle of advection, averaged over a set of initial conditions and CFL numbers $\nu \in [0.2, 0.8]$. The initial conditions include both smooth harmonic profiles from the test dataset and composite profiles combining discontinuities, sharp gradients, smooth regions, and extrema (see Appendix~\ref{sect:appendix-D}). The resulting score therefore reflects the ability of a limiter to preserve a broad range of shapes. Overall, the set has a substantial spread in performance, with several members achieving accuracy comparable to the classical OSTVD3 scheme (see Fig.~\ref{fig:ensemble_ddl} which shows all members in the Sweby diagram), while low-ranked members are excessively diffusive, damping extrema and smoothing sharp features of the solution. All learned members share common features, which are reflected in the best-performing limiter (Table~\ref{tab:ddl}) and discussed in the next section.

\begin{table}[]
    \centering
    \caption{Shape-preservation metric for the best data-driven limiter (DDL) from the set,  the reference OSTVD3 scheme \citep{arora_well-behaved_1997} and the Piecewise Parabolic Method (PPM) \citep{colella_piecewise_1984}, for comparison. Values correspond to the mean $L^1$ error after one cycle of advection averaged over CFL numbers $\nu \in [0.2, 0.8]$, for two classes of initial conditions: smooth harmonic profiles (testing dataset) and composite shapes combining Gaussian, square, triangular, and semi-elliptic features. The best one is bolded.}
    \begin{tabular}{llll}
        Limiter &  Testing Dataset & Composite shapes & Total\\
        DDL     & \textbf{7.42} $\times 10^{-2}$ & \textbf{1.17} $\times 10^{-2}$ & \textbf{4.29} $\times 10^{-2}$ \\
        OSTVD3  & 1.11 $\times 10^{-1}$ & 1.23 $\times 10^{-2}$ & 6.15 $\times 10^{-2}$ \\
        PPM     & 1.29 $\times 10^{-1}$ & 1.28 $\times 10^{-2}$ & 7.09 $\times 10^{-2}$ \\
    \end{tabular}
    \label{tab:ddl}
\end{table}

\subsection{Preservation of sharp gradients and extrema}\label{sect:sensitivity}

\begin{figure}
    \centering
    \includegraphics[width=0.8\linewidth]{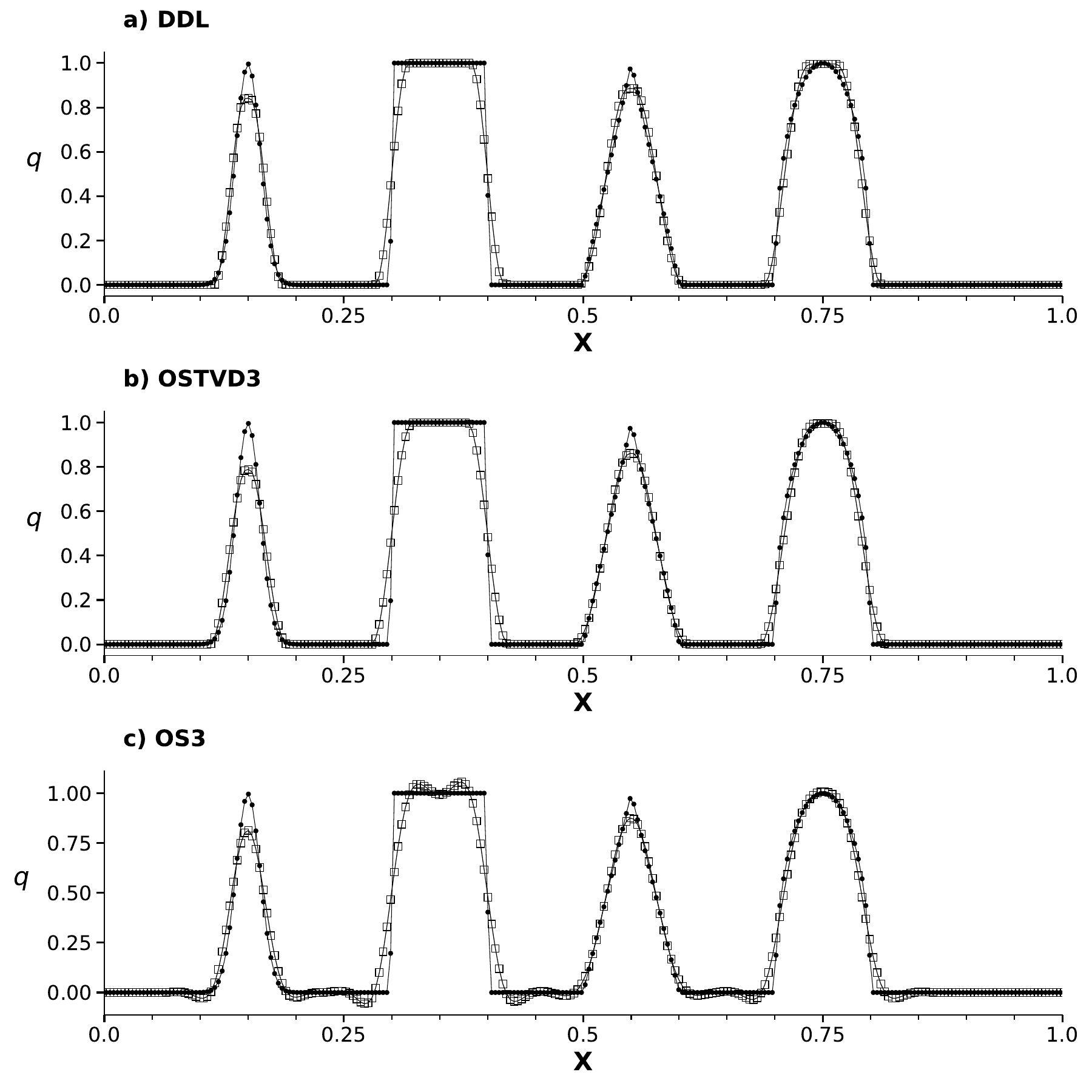}
    \caption{Advection over one period. The dotted line denotes the initial condition and the square markers denote the final state. From top to bottom: the optimized data-driven limiter (DDL), the TVD third order one-step scheme (OSTVD3) and the unlimited scheme (OS3), for comparison. All the schemes are used with CFL $\nu=0.2$.}
    \label{fig:ddl_signal1}
\end{figure}

\begin{figure}
    \centering
    \includegraphics[width=0.8\linewidth]{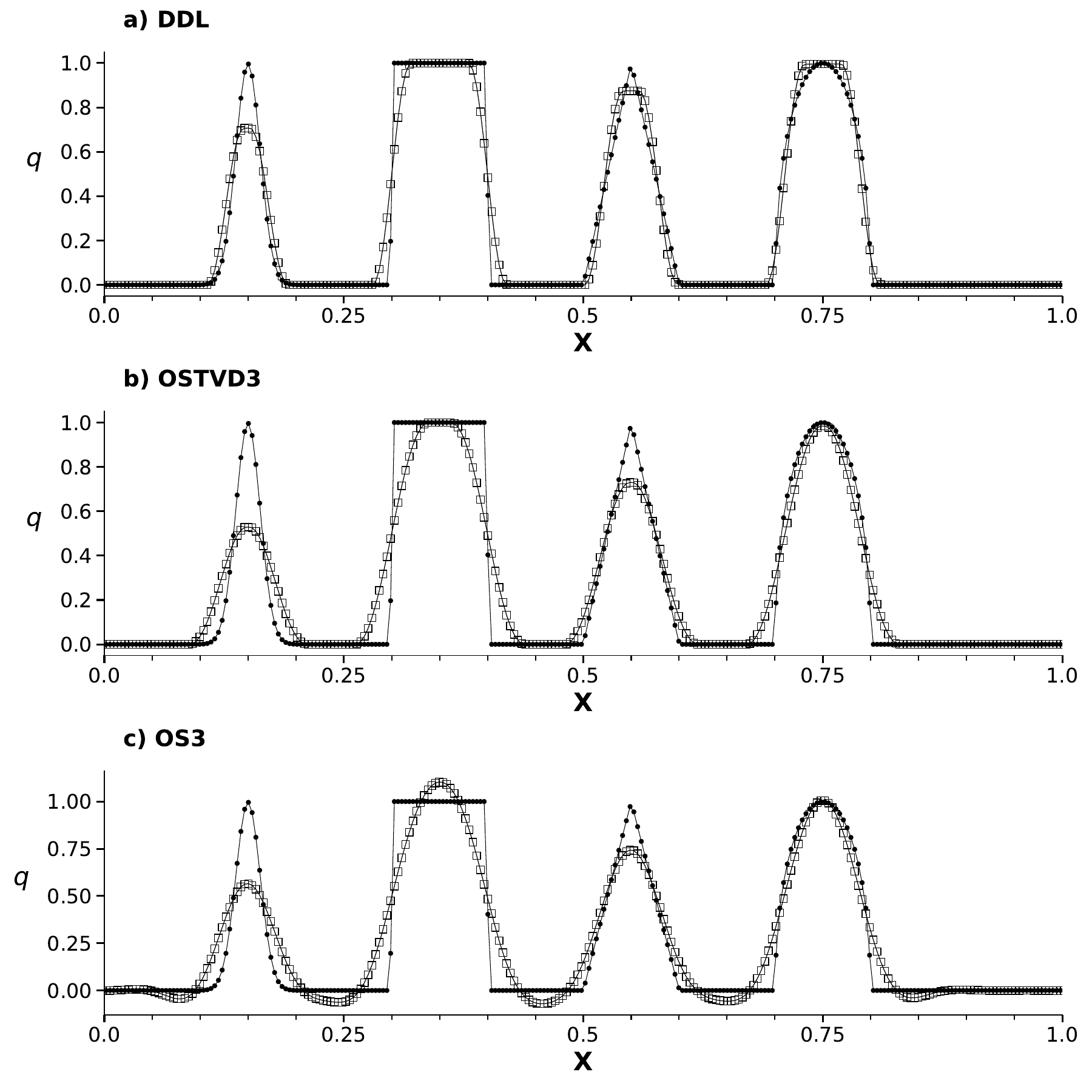}
    \caption{Advection over 10 periods. From top to bottom: the optimized data-driven limiter (DDL), the TVD third order one-step scheme (OSTVD3) and the unlimited scheme (OS3), for comparison. All the schemes are used with CFL $\nu=0.2$.}
    \label{fig:ddl_signal10}
\end{figure}

Fig.~\ref{fig:ddl_signal1}a illustrates the performance of the best data-driven limiter identified in Table~\ref{tab:ddl} (DDL) after one period of advection of the composite initial condition in Eq.~\ref{eq:dt}. The third-order unlimited scheme (OS3) performs well in smooth regions but generates oscillations near discontinuities (Fig.~\ref{fig:ddl_signal1}c). The TVD limiter removes these oscillations but diffuses extrema, although this effect is spatially localized (Fig.~\ref{fig:ddl_signal1}b). In contrast, the DDL preserves monotonicity while substantially reducing diffusion of extrema. Thus, the data-driven limiter yields a superior solution compared to the classical third-order OSTVD3 scheme. After ten periods (Fig.~\ref{fig:ddl_signal10}), these conclusions remain unchanged. A notable difference, however, is that extrema under the DDL generate plateaus more prominently than under TVD (compare Fig.~\ref{fig:ddl_signal10}a–b).

\subsection{Analysis of the dominant limiting mechanism}

\begin{figure}
    \centering
    \includegraphics[width=0.9\linewidth]{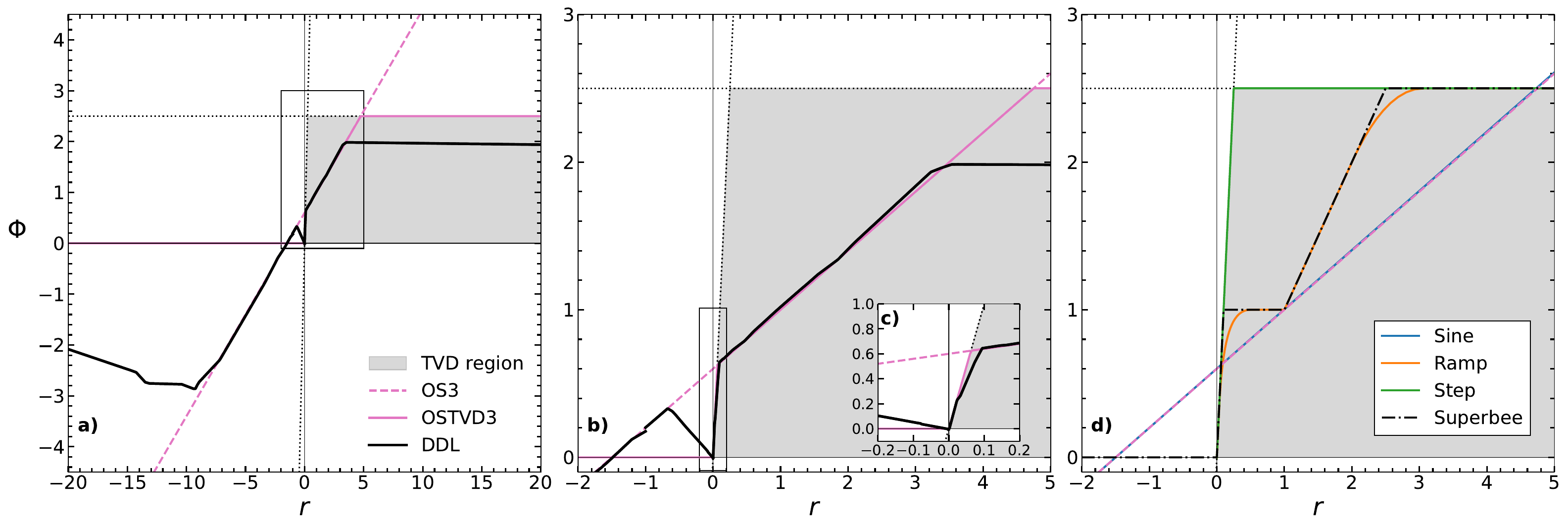}
    \caption{Data-driven and conventional flux limiters in the Sweby $r$-diagram. Panel (b) shows a magnified view of the inset in panel (a), while panel (c) further zooms into panel (b). The accuracy function $\Phi$ is plotted along the vertical axis and the slope ratio $r$ along the horizontal axis, with CFL $\nu=0.2$. The best-performing data-driven limiter DDL is shown as a solid black line. Panel d) shows the accuracy function, inverted from Eq.\ref{eq:cphi}, of a sine wave (blue), a linear ramp (orange), and an abrupt step (green). For reference, TVD region is shaded in gray and some conventional discrete schemes are shown in pink and black lines.}
    \label{fig:ddl}
\end{figure}

The observed compressive effect is more visible in the vicinity of smooth extrema, which may reflect that the neural network has learned a less diffusive reconstruction in non-monotonic regions. However, classical limiters such as Superbee also induce compression while reverting to first-order behavior near extrema \citep{leveque_numerical_1992}, raising the question of which underlying mechanism is responsible for the staircasing effect observed in DDL. To address this question, we analyze DDL by diagnosing the accuracy function $\Phi$ (obtained by inverting Eq.~\ref{eq:cphi}) and plot it as a function of the slope ratio $r$ in the Sweby diagram (Fig.~\ref{fig:ddl}a–c). The learned limiter aligns closely with the third-order OS3 scheme (dashed pink line) and remains within the TVD region (gray zone) for $r>0$. For $r<0$, it departs from classical TVD behavior by displaying a positive kink, followed by a transition toward a third-order-like behavior that ends in a negative plateau. Note that DDL exhibits a discontinuity near $r = -1$, due to approximation errors of the neural network at extreme curvature values $\delta^2 \tilde q = \pm 2$, however, it has no visible impact on the simulated solution.

\begin{figure}
    \centering
    \includegraphics[width=0.8\linewidth]{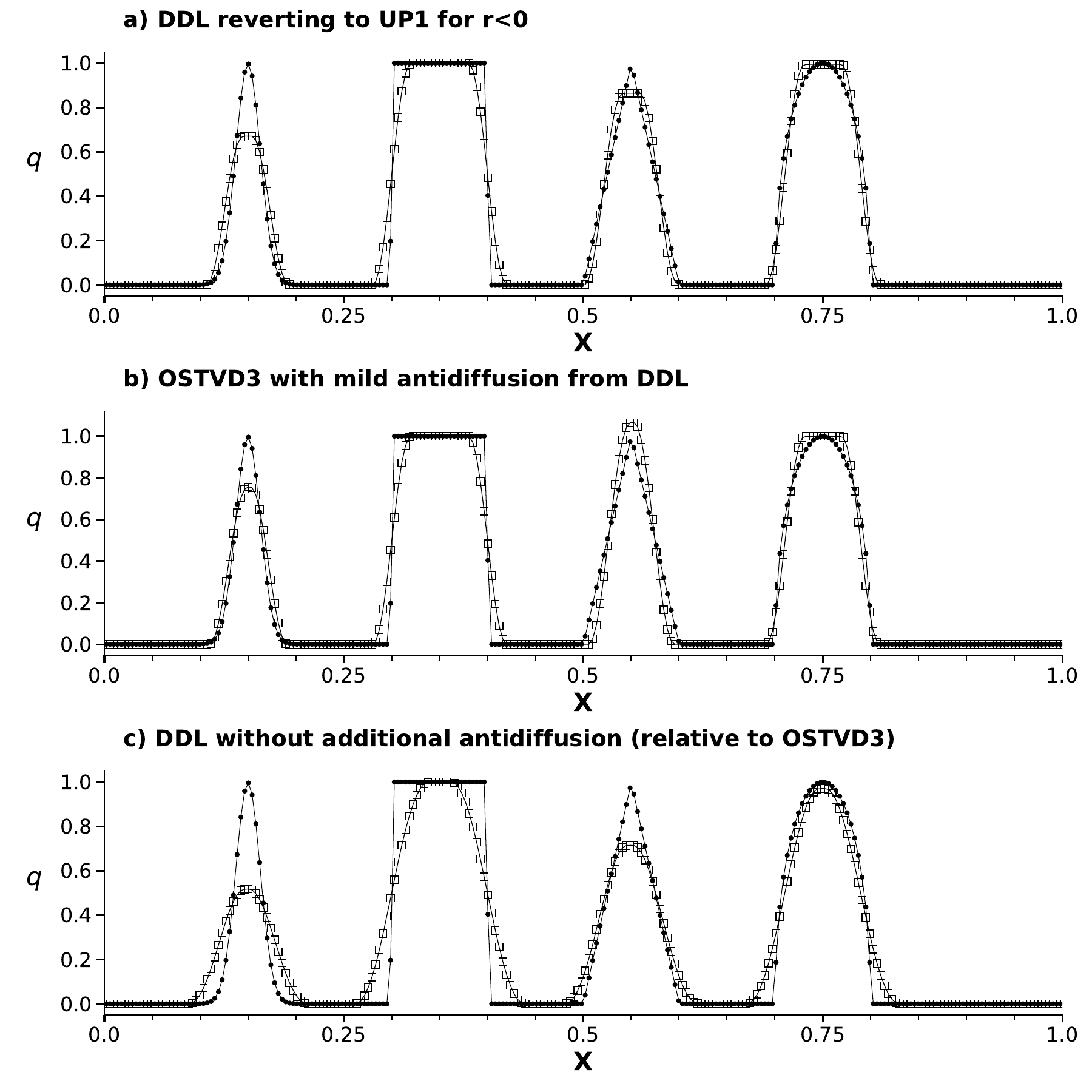}
    \caption{Advection over 10 periods for modifications of the DDL limiter. In (a), all $\Phi$ values for $r<0$ are set to zero. In panels (b-c), the limiter enforces DDL when $\Phi_{\text{DDL}}$ value is respectively above or below $\Phi_{\text{OSTVD3}}$ for $r>0$. All simulations use Courant number $\nu=0.2$.}
    \label{fig:ddl_sens_10}
\end{figure}

Fig.~\ref{fig:ddl_sens_10}a presents sensitivity test using modified versions of DDL where the entire $r<0$ branch is suppressed. After ten periods, the resulting solutions are nearly indistinguishable from the unmodified DDL (Fig.~\ref{fig:ddl_signal10}a), suggesting that the dominant mechanism lies in the $r\geq0$ regime. Next, comparing DDL (in black) with OSTVD3 (in pink) in Fig.~\ref{fig:ddl}b–c reveals small overshoots in near-linear regimes and undershoots in curved regimes, corresponding to mild antidiffusion and additional diffusion relative to OSTVD3, respectively. Sensitivity tests show that retaining only the antidiffusive excess relative to OSTVD3 is sufficient to reproduce the observed staircasing effect (Fig.~\ref{fig:ddl_sens_10}b), whereas additional diffusion is required to prevent out-of-range extrema. In summary, the staircasing behavior in DDL is due to mild antidiffusion in near-linear curvature regimes, which significantly improves the preservation of extrema and sharp fronts. 

\subsection{Numerical convergence}

While mild antidiffusion is optimal for DDL, it remains unclear whether the benefit of better preserving discontinuities and steep gradients is detrimental to the accuracy in smooth regions. To assess this trade-off, Table~\ref{tab:accddl} reports the convergence order of the limited schemes for a smooth initial condition. The theoretical third-order scheme with TVD constraints (OSTVD3) converges at nearly second order (1.96), whereas the data-driven limiter DDL achieves mean accuracy orders of 1.78. This reduction in convergence rate is consistent with the observed staircasing effect, which is most pronounced for smooth solutions. However, more work is needed to better understand the impact of flux limiting on the convergence rate.

\begin{table}[]
    \centering
    \caption{Transport of a sum of 4 sine waves: $L_2$ error and order of accuracy for the data-driven flux-limiter DDL and the one-step third order scheme combined with TVD constraints (OSTVD3).}
    \begin{tabular}{llll}
         Method & Number of grid points & $L_2$ error & $L_2$ order \\
         DDL &  160  & 2.689 $\times 10^{-2}$ &     \\
    &  320  & 7.768 $\times 10^{-3}$ & 1.79 \\
    &  640  & 2.310 $\times 10^{-3}$ & 1.75 \\
    & 1280  & 6.148 $\times 10^{-4}$ & 1.91 \\
    & 2560  & 1.950 $\times 10^{-4}$ & 1.66 \\ 
    &&&\\
            OSTVD3 &  160  & 3.029 $\times 10^{-2}$ &     \\
    &  320  & 8.076 $\times 10^{-3}$ & 1.91 \\
    &  640  & 2.318 $\times 10^{-3}$ & 1.80 \\
    & 1280  & 5.060 $\times 10^{-4}$ & 2.20 \\
    & 2560  & 1.238 $\times 10^{-4}$ & 2.03 \\ 
    \end{tabular}
    \label{tab:accddl}
\end{table}

\subsection{Interpretation of distortion effect for traditional limiters}

The numerical flux associated with each profile: a discontinuous step, a linear ramp, and a smooth sine wave can be expressed through the corresponding accuracy function $\Phi$ (obtained by inverting Eq.~\ref{eq:cphi}), which is shown in Fig.~\ref{fig:ddl}d. First, the compressive Superbee limiter closely matches the ramp profile (in orange), consistent with its design based on the method of characteristic for linear jumps \cite{roe_characteristic-based_1986}. Further, the step profile (in green) delineates the boundary of the TVD region (shaded gray), indicating that Harten’s TVD bounds in Eq.~\ref{eq:harten} correspond to the limiting behavior associated with advecting sharp fronts. Lastly, we note that the accuracy functions associated with the sine wave (in blue) and the step intersect, suggesting that the emergence of spurious oscillations seen in Fig.~\ref{fig:learned}d and Fig.~\ref{fig:ddl_signal1}c arises from the polynomial reconstruction of a discontinuous front based on discrete inputs, which are treated as if it was those of a smooth solution.

\section{Synthesis and Discussion}

\subsection{Synthesis}
\begin{figure}
    \centering
    \includegraphics[width=0.5\linewidth]{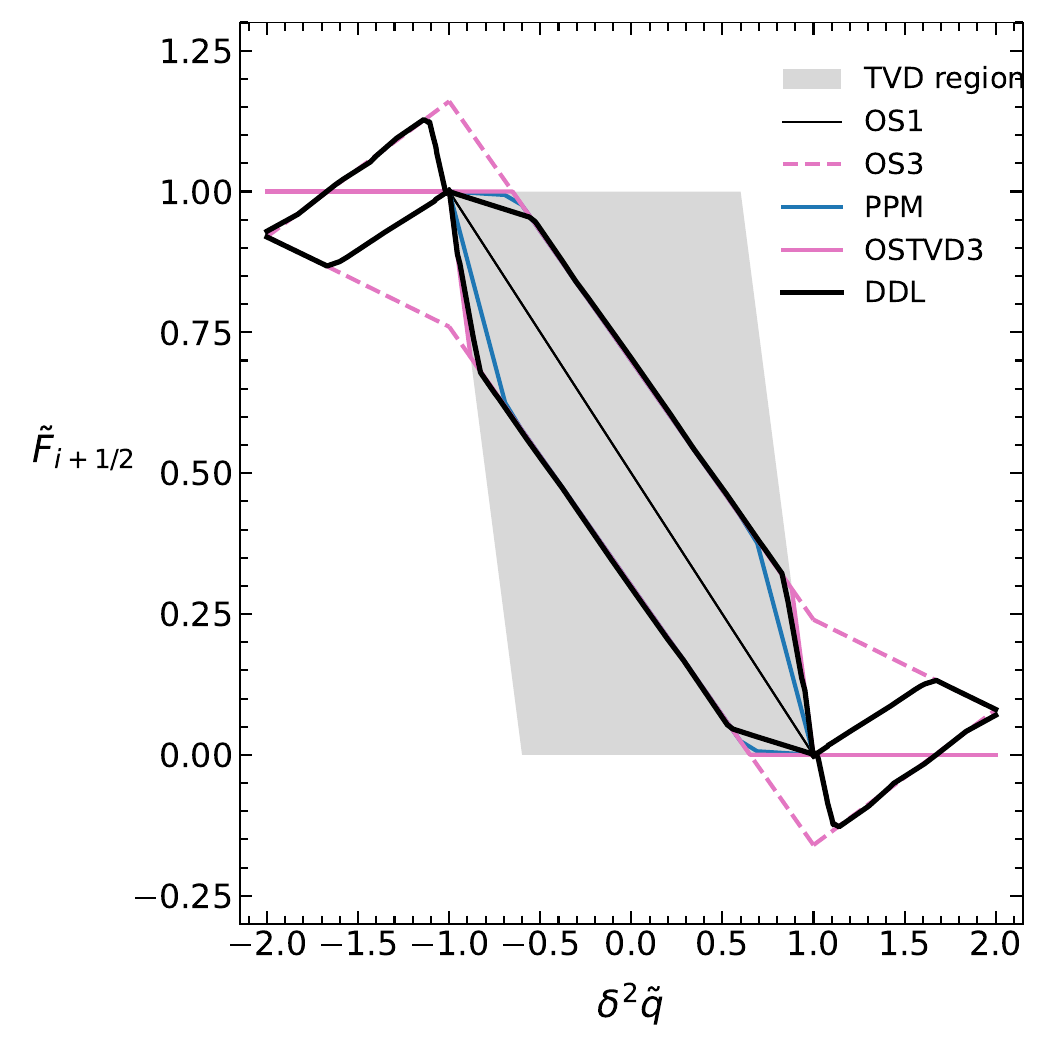}
    \caption{Advective schemes diagnosed in the flux-curvature diagram. Each color denotes a shape specified in the legend. The normalized numerical flux $\tilde F_{i+1/2}^n$ is plotted in the $y$-axis, for a CFL number $\nu=0.2$, against the local curvature of the normalized gridded field, $\tilde q_{i-1}-2\tilde q_{i}+\tilde q_{i+1}$ along $x$-axis.}
    \label{fig:synthesis}
\end{figure}

The resulting data-driven limiter DDL is summarized in Fig.~\ref{fig:synthesis}, where the predicted flux is plotted in a flux-curvature diagram. For comparison, OSTVD3 and the Piecewise Parabolic Method (PPM) \citep{colella_piecewise_1984} are also shown. Notably, DDL (thick black line) remains close to the analytical OSTVD3 (in pink) and PPM (in blue) curves for monotonic region ($r\geq 0$), highlighting how small modifications to the limiter design can significantly improve shape-preservation. Both OSTVD3 and PPM naturally revert to the first-order upwind scheme (thin black line) in non-monotonic regions, visible as plateaus in the diagram. Thus, representing the numerical flux as a function of the curvature of the discrete cell-averages provides a simple visualization of local reconstructions using a three-point stencil. In this framework, both slope limiters and flux limiters correspond to specific trajectories in the diagram, where slope breaks are consistent with limiting choices. In short, Fig.~\ref{fig:synthesis} shows that limiter design can be interpreted, and potentially guided, through the geometry of these local reconstructions.

\subsection{Discussion}
All randomly seeded $F_\text{NN}$ reliably had a non-zero feature in the non-monotonic region of $r<0$ (Fig.~\ref{fig:ensemble_ddl}) and yet, when eliminated during inference, this region had no appreciable effect. We are unable to explain why this feature is selected in training when it is not substantially affecting the solution. The DDL slopes match the OS3 reconstruction in some cases for non-monotonic discrete inputs (e.g in Fig.~\ref{fig:synthesis}), echoing \citep{woodfield_new_2024,spekreijse_multigrid_1987} who identified admissible reconstructions for non-monotonic inputs that satisfy TVD constraints while differing from first-order upwind scheme. However, the sensitivity analysis in Section~\ref{sect:sensitivity} suggests that allowing higher-than-first-order reconstruction for $r<0$ has limited impact  (compare extrema heights in Fig.~\ref{fig:ddl_sens_10}a with Fig.\ref{fig:ddl_signal10}a). We suspect that improvements may emerge when considering larger stencils, which could better resolve non-monotonic regions and mitigate the known limitations of TVD limiters near critical points \citep{daru_high_2004,colella_limiter_2008}.

While this study identifies a data-driven flux limiter that would be difficult to derive through analysis alone, it does not preclude the existence of alternative optimal limiters under different design criteria. In particular, extending the learning framework to incorporate metrics beyond point-wise accuracy, such as monotonicity or entropy consistent constraints \citep{patel_thermodynamically_2022,romemont_data-driven_2025}, offers promising avenues for further exploration. More broadly, our results show that robust data-driven discretizations can be constructed when the underlying formulation is grounded in a stable numerical framework. Here, stability and generalization are ensured by a one-step formulation \citep{lax_systems_1960} and by stencil-scale normalization, provided that the neural network is adequately trained.

In this approach, the neural network infers the numerical flux by taking as inputs the local curvature $\delta^2 \tilde q$ together with a mirroring condition. This contrasts with conventional analyses, where slope ratios are used directly as parameters for constructing high-order and more elaborate limiters \citep{daru_high_2004, del_pino_arbitrary_2006, colella_limiter_2008}. In particular, the slope ratio $r$ is invariant under local affine transformations, ensuring the semilinearity of $\Phi$, which is a necessary condition for the stability and consistency of the discretization. However, using an unbounded quantity such as $r$ as input to a neural network can be problematic from a machine-learning perspective, as it may induce vanishing gradients or out-of-distribution behavior \citep{zhuang_learned_2021}. Sensitivity tests with the representation $\Phi(r,\nu)$ gave similar results and confirmed that large values of $r$, associated with sharp gradients, occasionally induce numerical instabilities, even when clipping strategies were applied (not shown). This likely reflects the inability of the training dataset to guarantee robust behavior outside the sampled range. Last, an alternative is to use as inputs the local cell-averaged values $(q_i)_{i\in K}$ combined with a mirroring condition. However, this increases the input dimensionality by two, which typically entails larger training datasets and network capacity to achieve comparable accuracy \cite{goodfellow_deep_2016}. In summary, finding which input representation provides the best performance remains a critical step when designing a data-driven discretization, which should be assessed based on the application.

Finally, numerical simulations always involve a trade-off between accuracy and computational cost. Although, neural-network inference relies primarily on dense matrix multiplications and therefore typically requires more operations than traditional discretizations, practical performance depends strongly on hardware characteristics, including memory access patterns and processor architecture, which can significantly affect runtime \cite{bar-sinai_learning_2019, zhuang_learned_2021}. Naturally, it is not computationally efficient to preserve at all cost the full neural network. Hence, to mitigate computational overhead, we perform a model-order reduction of DDL by constructing a parametric representation of the flux limiter based on hinge functions with CFL-dependent knot locations and polynomial amplitudes (presented in Appendix~\ref{sect:appendix-E}), following classical approaches in regression splines \cite{hastie_elements_2009,friedman_multivariate_1991}. The approximated limiter reproduces well the results obtained with the full neural-network DDL (not shown). The ability to reduce the model complexity suggests that we might one day be able to connect data-driven improvements to numerical analysis and theory, which we did not do here.

\section{Conclusions}

We investigated design principles for stable and generalizable data-driven parameterizations of the linear advection problem. Within a finite-volume framework for one-dimensional transport in a periodic domain, a neural network was constructed from first principles and used as a numerical advection scheme. Several data-driven strategies were examined for discrete advection and flux limiting across a range of initial conditions, from discontinuous to smooth profiles.

(i) We show that exact advection of profiles spanning discontinuous to smooth shapes requires distinct discretizations when only cell-averaged values are used as inputs. From a machine-learning perspective, mixing profiles with widely varying curvature creates a multi-valued learning problem, which explains the limited generalization of data-driven schemes that mimic conventional finite-volume reconstructions. For smooth data, where curvature is well resolved at a given Courant number, the associated numerical flux is already well approximated by classical high-order schemes. More broadly, it is always possible to derive a discretization that exactly fits a given solution shape, implying that data-driven approaches are unlikely to outperform analytical discretizations so long as the loss function reflects accuracy in the same way as usually defined in numerical analysis. 

(ii) It is shown that numerical stability can be achieved through a rigorous formulation of the numerical framework and by enforcing the semilinearity property of the numerical-flux, ensuring invariance under affine transformations of the inputs. While previous studies enforced semilinearity through global normalization, we demonstrate that applying normalization locally at the stencil scale leads to more stable and generalizable data-driven discretizations. However, the choice of training data remains the dominant factor in achieving stable and accurate discretizations.

(iii) We introduce a new data-driven flux limiter, expressed as functions of the local curvature of the three-point stencil $\delta^2 q$ and the Courant number, which outperforms the classical OSTVD3 limiter in terms of shape preservation and generalizes well across smooth and discontinuous initial conditions. The results indicate that limiter design is particularly sensitive in near-linear curvature regimes, where mild antidiffusion significantly improves the preservation of extrema and sharp gradients. By contrast, allowing higher-order reconstruction in non-monotonic regions rather than reverting to first-order upwinding does not yield a clear improvement for overall advection accuracy.

(iv) We show that training on polynomial profiles defined over the stencil yields stable, high-order accurate data-driven discretizations, with the polynomial degree directly controlling the formal order of accuracy. Achieving sixth-order accuracy or higher, however, requires careful control of neural-network approximation errors. These results highlight promising pathways for constructing multidimensional data-driven advection schemes that are both accurate and stable.

We conclude that finite-volume advection schemes based on data-driven reconstructions using only cell-averaged inputs tend to specialize toward the solution shapes represented in the training dataset. However, conventional schemes can also exhibit implicit specialization toward particular solution classes. Addressing the generalization limits of data-driven discretizations therefore connects to the broader challenge of developing robust advection schemes in numerical models and may provide new directions for improving the fidelity of simulated currents across diverse flow regimes.

%
%

\section*{Acknowledgments}
The authors thank Joseph Mouallem for comments on the draft of this manuscript. The authors would like to also thank M$^\text{2}$LInES team in particular Pavel Perezhogin, Danni Du and Will Chapman, for fruitful discussion and comments. This research received support through Schmidt Sciences, LLC, under the M$^\text{2}$LInES project. The code for training of the neural networks and figures is available at \url{https://github.com/antoine-182/design-principles-ML-FV-schemes-JCP2026}.


\appendix

\renewcommand{\theequation}{\Alph{section}\arabic{equation}}
\renewcommand{\thefigure}{\Alph{section}\arabic{figure}}
\renewcommand{\thetable}{\Alph{section}\arabic{table}}

\section{Computation of the normalized FV flux for analytical profile shapes}
\label{sect:appendix-A}
\setcounter{table}{0}    
\setcounter{figure}{0}
\setcounter{equation}{0}

\subsection{Principle} \label{sect:appendix-A1}
Here, we present how we compute the learnable flux (of section \ref{sect:multivalued}) shown in the flux-curvature diagram Fig.~\ref{fig:learnable}. We wanted to gain insights about the advection of some well-know shapes from the perspective of a discrete scheme, to understand better how to learn this with a neural network. This means that we have to compute the FV flux in function of the cell-averaged values in every point of the solution. Both finite-volume cell-averages values and numerical flux are consistent to profile shape, which, if learned by a neural network, enables perfect advection of the same solution. 

i) To achieve this, we computed the pairs of flux and cell-averaged solution by sampling the profile shape, that is shifting incrementally the full profile across the length of the stencil. The logic of the sampling method is reflected in the following simple example of a unit positive jump at $x=w_0$:
\begin{eqnarray}
  f_{jump}(x) =
\begin{cases}
1 & x \ge w_0, \\
0 & \text{otherwise.}
\end{cases}  
\end{eqnarray}
Sampling consists of calculating the FV flux in Eq.~\ref{eq:cflux} and the cell-averages with Eq.~\ref{eq:cavg} for every increment along the $x$-dimension. The general expression for the numerical flux, defined in $x_{i+1/2}$, is
\begin{eqnarray}
  F_{jump}(x;w_0) =
\begin{cases}
1 & x - w_0 \ge \nu \Delta, \\
\frac{x-w_0}{\nu \Delta} &   0 \le x-w_0 \leq \nu \Delta, \\
0 & \text{otherwise.}
\end{cases}  
\end{eqnarray}
and the cell-averaged solution, defined in $x_{i}$, is 
\begin{eqnarray} \label{eq:cavg_jump}
    q_{jump}(x; w_0) = \text{minmod}\big(\frac{1}{\Delta}(x-w_0) - \frac{1}{2} \big)
\end{eqnarray}

ii) Considering  a three-node scheme, which is the minimal stencil size for capturing curved inputs, we concatenate the formulation Eq.~\ref{eq:cavg_jump} to form $(q_{i-1},q_i,q_{i+1})$. To obtain the invariant form of the distribution, both flux and the input cell-averaged stencil are normalized with the min-max rule (see Section~\ref{sect:minmax}). 

iii) Using the fact that each permutations $(\tilde q_{i-1},\tilde q_{i},\tilde q_{i+1})$ can be described by the local curvature of the gridded data, $\delta^2_{i+1/2}\tilde q = \tilde q_{i-1} - 2\tilde q_{i} + \tilde q_{i+1}$, and the sign of the slope $\tilde q_{i+1}-\tilde q_{i-1}$ (see Fig.~\ref{fig:K3}), we partition the FV flux as a function of these two parameters. The case of uniform values with the normalization reduces to the trivial case which leads to a singular point in the phase space and thus is ignored in the plot.

\subsection{Training/test profiles considered in the study}

While previous efforts have used numerical quadrature to estimate the correct flux, to be able to test numerical agreement we calculate the exact flux analytically. We present the exact expression of the FV flux and the cell-averaged solution for three analytical profiles shown in Fig.~\ref{fig:learnable} commonly used for training: a discontinuous step, a linear ramp, and a sine wave. 

\subsubsection*{Discontinuous Step} 
A discontinuous step defined by $w_1>w_0$, large enough so that the width of the plateau is wider than $K\Delta$, is
\begin{eqnarray} \label{eq:square}
  f_{step}(x) =
\begin{cases}
1 &  w_0 \leq x \leq w_1,  \\
0 & \text{otherwise,}
\end{cases}  
\end{eqnarray}
which gives the flux
\begin{eqnarray}
F_{step}(x;w_0,w_1) =
\begin{cases}
1  & w_0 + \nu\Delta \leq x \leq w_1, \\
1 - \frac{x-w_1}{\nu \Delta} &  0 \le x-w_1 \leq \nu \Delta, \\
\frac{x-w_0}{\nu \Delta} &  0 \le x-w_0 \leq \nu \Delta, \\
0 & \text{otherwise,}
\end{cases}   
\end{eqnarray}
and cell-averaged solution
\begin{eqnarray}
    q_{step}(x;w_0,w_1) = \text{minmod}\big(\frac{1}{\Delta}(x-w_0) - \frac{1}{2} \big) - \text{minmod}\big(\frac{1}{\Delta}(x-w_1) - \frac{1}{2} \big).
\end{eqnarray}

\subsubsection*{Linear Ramp}  
We consider a linear ramp of slope larger than the size of a cell of the mesh $\Delta$. Expressions are given only for the increasing portion of the profile, the negative portion can be obtained in a similar way. The profile shape is 
\begin{eqnarray}
  f_{ramp}(x) =
\begin{cases}
1 &   x \geq w_1,  \\
\frac{x}{w_1 - w_0}  &   w_0 \geq x \geq w_1  \\
0 & \text{otherwise}
\end{cases}  
\end{eqnarray}
which gives a quadratic dependency with the position of the ramp. The numerical flux is
\begin{eqnarray}
  F_{ramp}(x;w_0,w_1) =
\begin{cases}
1  &  x \geq w_1 + \nu \Delta, \\
\frac{x-w_1}{\nu\Delta}+\frac{w_1^2 - (x-\nu\Delta)^2}{2\nu \Delta(w_1-w_0)} &  x \leq w_1 + \nu \Delta, \quad\text{and}\quad x \geq w_1 , \\
\frac{2x-\nu \Delta}{2(w_1-w_0)} &  x \geq w_0 + \nu \Delta, \quad\text{and}\quad x \leq w_1, \\
\frac{x^2-w_0^2}{2\nu \Delta(w_1-w_0)} &  x \leq w_0 + \nu \Delta, \quad\text{and}\quad x \geq w_0, \\
0  & \text{otherwise,}
\end{cases}  
\end{eqnarray}
and cell-averaged solution is
\begin{eqnarray}
q_{ramp}(x;w_0,w_1) =
\begin{cases}
1  &  x \geq w_1 +  \Delta/2, \\
\frac{1}{2} + \frac{x-w_1}{\Delta} +\frac{w_1^2 - (x- \Delta/2)^2}{2 \Delta(w_1-w_0)} &  x \leq w_1 +\Delta/2 , \quad\text{and}\quad x \geq w_1 -\Delta/2, \\
\frac{x}{ (w_1-w_0)} &  x \geq w_0 + \Delta/2, \quad\text{and}\quad x \leq w_1-\Delta/2, \\
\frac{(x+\Delta/2)^2-w_0^2}{2\Delta(w_1-w_0)} &  x \leq w_0 + \Delta/2, \quad\text{and}\quad x \geq w_0 - \Delta/2, \\
0  & \text{otherwise.}
\end{cases}  
\end{eqnarray}

\subsubsection*{Sine wave} 
A sine wave, of wavelength much larger than the Nyquist mode is
\begin{eqnarray} \label{eq:sine}
  f_{sine}(x) = \sin{(k x)}
\end{eqnarray}
with wavenumber $k=2\pi/\lambda$. Following Eq.~\ref{eq:cflux} and Eq.~\ref{eq:cavg}, the finite-volume flux is 
\begin{eqnarray}
  F_{sine}(x) = \sinc{\left(\frac{k \nu \Delta}{2}\right)}\sin{\left(kx-\frac{k\nu \Delta}{2}\right)}
\end{eqnarray}
and the cell-centered averaged solution writes
\begin{eqnarray}
    q_{sine}(x)= \sinc{\left(\frac{k \Delta}{2}\right)}\sin{\left(kx-\frac{k\Delta}{2}\right)}
\end{eqnarray}

\subsubsection*{Normalization step}
We normalize both the flux and the stencil data to obtain the invariant form of the discretization. The cases of the discontinuous step and the linear ramp are unitary, therefore the expressions remain unchanged after min-max normalization. When considering a general profile, we can follow a slightly different approach for the sampling: the maxima $(q_{\min},q_{\max})$ across the stencil are derived depending on the monotonicity of the solution, then the normalized flux is obtained by $\widetilde F = (F- q_{\min})/ (q_{\max}-q_{\min})$. Hence, for the sine wave, $f_{sine}$, the specific normalization parameters are
\begin{eqnarray}
  (q_{\min},q_{\max})(x_i) =
\begin{cases}
\begin{aligned}
& (q_{i-1},q_{i+1}) & \forall \ x_i & \in [\frac{\Delta}{2}-\frac{ \lambda}{4};\frac{ \lambda}{4}-\frac{\Delta}{2}], \\
& (q_{i+1},q_{i-1}) & \forall \ x_i & \in [ \frac{ \lambda}{4}+\frac{\Delta}{2};\frac{3\lambda}{4}-\frac{\Delta}{2}], \\
& (q_{i-1},q_{i})   & \forall \ x_i & \in [ \frac{ \lambda}{4}-\frac{\Delta}{2};\frac{ \lambda}{4}                 ], \\
& (q_{i+1},q_{i})   & \forall \ x_i & \in [ \frac{ \lambda}{4}                 ;\frac{ \lambda}{4}+\frac{\Delta}{2}], \\
& (q_{i},q_{i-1})   & \forall \ x_i & \in [ \frac{3\lambda}{4}-\frac{\Delta}{2};\frac{3\lambda}{4}                 ], \\
& (q_{i},q_{i+1})   & \forall \ x_i & \in [ \frac{3\lambda}{4};\frac{3\lambda}{4}+\frac{\Delta}{2}],   
\end{aligned} 
\end{cases}
\end{eqnarray}
which gives, for the monotonic increasing case ($\delta^2 q \in [-1,1]$, blue line in Fig.~\ref{fig:learnable}a), the following analytical expression:
\begin{eqnarray}
  \forall\ x\in[\frac{\Delta}{2}-\frac{ \lambda}{4};\frac{ \lambda}{4}-\frac{\Delta}{2}],\quad
  \tilde F_{sin}(x) = \frac{\sinc\left(\frac{k \nu \Delta}{2}\right) \sin{\left( kx - \frac{k\nu \Delta}{2}  \right)}
                           -\sinc\left(\frac{k        \Delta}{2}\right) \sin{\left( kx - \frac{k 3\Delta}{2}  \right)} }
                           { 2 \sinc\left(\frac{k\Delta}{2}\right) \sin{(k\Delta)} \cos{\left( kx - \frac{k \Delta}{2}  \right)} }
\end{eqnarray}

\section{Hyper-parameters} \label{sect:appendix-B}
\setcounter{table}{0}    
\setcounter{figure}{0}
\setcounter{equation}{0}

The data-driven discretizations are modeled by multi-layer perceptrons (MLPs), i.e. fully connected neural networks. The architectures employed in each part of the study, together with their associated hyperparameters, are summarized in Tables~\ref{tab:hyperA}–\ref{tab:hyperB}–\ref{tab:hyperC}. All computations are performed in Python, and the data-driven discretizations are implemented using the PyTorch framework \citep{paszke_pytorch_2019}. Network parameters are optimized with the Adam optimizer using a learning rate of $10^{-3}$, and all training and inference are carried out in double-precision floating point. In all cases, the generated dataset is split into training and testing subsets using an 80/20 ratio, and training is performed with a batch size of 32 samples.

\begin{table}[]
    \centering
    \caption{Parameters used for the shape-specific neural-network presented in Section~\ref{sect:2.3}.}
    \begin{tabular}{ll}
         Inputs & $(q_{i-1}^n, q_{i}^n, q_{i+1}^n)$\\
         Output & $F_{i+1/2}^n$ \\
         Number of hidden layers& 3\\
         Number of filters& 8 \\
         Activation function & Relu \\
         Number of epochs & 300 \\ 
    \end{tabular}
    \label{tab:hyperA}
\end{table}

\begin{table}[]
    \centering
    \caption{Parameters used for the one-step data-driven discretization trained with polynomials presented in Section~\ref{sect:2.4}.}
    \begin{tabular}{llll}
                        & OSDD3 & OSDD4 & OSDD5 \\
    Inputs              & $(q_{i-1}^n, q_{i}^n, q_{i+1}^n)$ & $(q_{i-2}^n,q_{i-1}^n, q_{i}^n, q_{i+1}^n)$ & $(q_{i-2}^n,q_{i-1}^n, q_{i}^n, q_{i+1}^n,q_{i+2}^n)$ \\
    Output              & $F_{i+1/2}^n$ & $F_{i+1/2}^n$ & $F_{i+1/2}^n$ \\
    Number of hidden layers& 1& 1& 1\\
    Number of filters   & 8 & 8 & 16\\
    Activation function & Linear & Linear & Linear \\
    Number of epochs    & 200 & 200 & 400 \\ 
    \end{tabular}
    \label{tab:hyperB}
\end{table}

\begin{table}[]
    \centering
    \caption{Parameters used for the data-driven limiter approach presented in Section~\ref{sect:ddfl}. The network predicts a correction to the first-order upstream (UP1) flux to improve accuracy.} 
    \begin{tabular}{l|l}
    Inputs              & $(\delta^2 q,\nu)$\\
    Output              & $\Delta F_{i+1/2}^n$\\
    Number of hidden layers& 3\\
    Number of filters& 32 \\
    Activation function & Relu \\
    Number of epochs    & 40\\ 
    \end{tabular}
    \label{tab:hyperC}
\end{table}

\section{Dataset} \label{sect:appendix-C}
\setcounter{table}{0}    
\setcounter{figure}{0}
\setcounter{equation}{0}

\subsection*{Shape-specific neural networks used in Section~\ref{sect:2.3}}
Neural networks shown in panel b,c and d of Fig.~\ref{fig:learned} are trained against a step profile of width $w_1-w_0 = 11\Delta$ (Eq.~\ref{eq:square}) and a sine wave of wavelength $L/8$ (Eq.~\ref{eq:sine}), respectively. The dataset consists of 4000 pairs of inputs/outputs $((q_{i-1},q_i,q_{i+1}); F_{i+1/2})$ by sampling the analytical functions using numerical interpolation, see Appendix~\ref{sect:appendix-A}.
 
\subsection*{One-step data-driven (OSDD) schemes trained with polynomials}
One-step data-driven discretizations achieves formal accuracy by training on polynomials defined on the stencil lengths. For training OSDD schemes, 2000 local polynomials of maximum degree $k-1$ are generated with random coefficients. Let $P_k$ denote a polynomial of degree $k\in\mathbb{N}^+$,
\begin{eqnarray}
P_k(x) = \sum_{l=0}^{k} \beta_l x^l,
\end{eqnarray}
where $\beta_l$ are scalar coefficients and $x\in\mathbb{R}$ is the spatial coordinate. Any affine transformation of the spatial coordinate, $x\mapsto \theta_1 x+\theta_2$, maps $P_k$ onto another polynomial of the same degree. As a result, the polynomial generation is independent of the underlying spatial grid and of the local grid spacing. Without loss of generality, we therefore set a uniform grid of unitary spacing centered at the midpoint of the first cell upstream of the flux. Following Eq.~\ref{eq:cavg} and~\ref{eq:cflux}, the cell-averaged concentration $q_i$ and the right interface mass flux $F_{1/2}$ are analytically calculated, thereby reducing numerical errors. They are given, respectively, by
\begin{eqnarray} 
q_i = \sum^k_{l=0} \frac{\beta_l}{l+1} \big((x_{i+1/2})^{l+1}-(x_{i-1/2})^{l+1}\big) \quad\text{ and }\quad F_{1/2} = \frac{u}{\nu}\sum^k_{l=0} \frac{\beta_l}{(l+1)} \big((x_{i+1/2})^{l+1}-(x_{i+1/2}-\nu)^{l+1}\big). 
\end{eqnarray}

\subsection*{Data-driven limiters}
To generate training data for section \ref{sect:ddfl}, we construct 53,248 harmonic solutions over a periodic domain of length $L = 2\pi$, with resolution of 256 grid points. Each realization is defined as the superposition of four sine modes with randomly sampled amplitudes and wavelengths, the latter constrained to be larger than 16 grid points to ensure smooth curves. The reference “true” solutions are generated by advecting each harmonic solutions on one timestep, and subsequently computing the next time step solutions using Eq.~\ref{eq:cavg}. The Courant number is sampled uniformly within the interval $[0.01,0.99]$ using 26 discrete values. 

\section{Sensitivity to randomized weight initialization} \label{sect:appendix-D}
\setcounter{table}{0}    
\setcounter{figure}{0}
\setcounter{equation}{0}

\begin{figure}
    \centering    
    \includegraphics[width=0.8\linewidth]{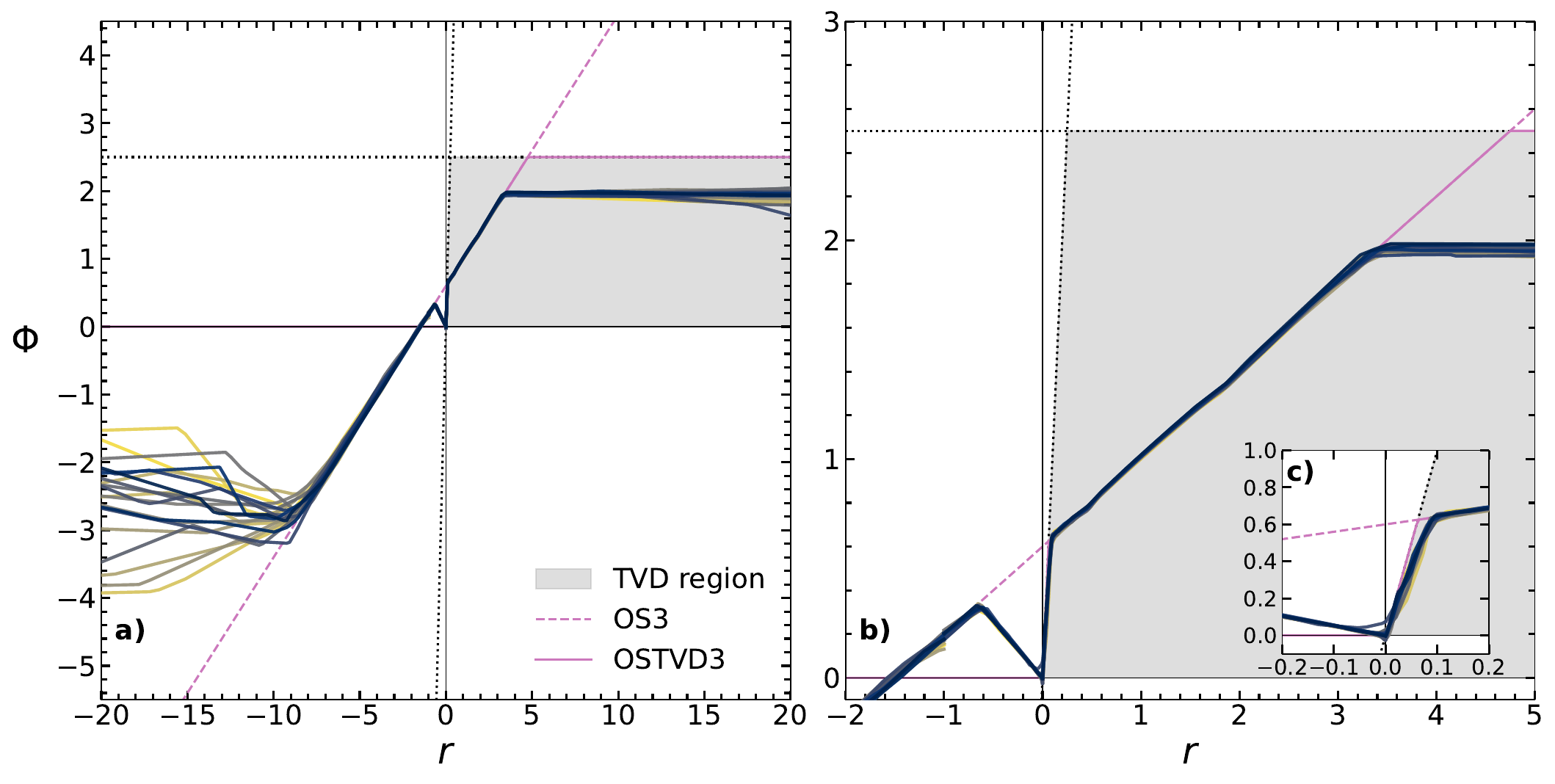}
    \caption{All members of the DDL formulation in the Sweby diagram plotted for a Courant number of 0.2. Colors range from blue (best) to yellow (worst), indicating the ranking based on shape-preservation performance.}
    \label{fig:ensemble_ddl}
\end{figure}

We quantify shape preservation by the mean $L^1$ error after one full advection cycle, averaged over CFL numbers $\nu \in \{0.2, 0.3, 0.4, 0.5, 0.6, 0.7, 0.8\}$ and a set of random initial conditions. The test suite consists of 16 samples drawn (i) from the test dataset and (ii) from a composite profile comprising a Gaussian pulse, a square wave, a triangular wave, and an ellipse, randomly shifted within the domain. The latter is a worthy testcase that combines discontinuities, sharp gradients, smooth regions, and extrema, and is defined on $[-1,1]$ using the transformed variable $x\mapsto(x/L-1/2)$ as
\begin{eqnarray} \label{eq:dt}
  q_0(x) =
\begin{cases}
\exp(-\log(2)(x+0.7)^2/0.0009) &  -0.8 \leq x \leq -0.6, \\
1                              &  -0.4 \leq x \leq -0.2, \\
1 - |10(x-0.1)|                &     0 \leq x \leq 0.2, \\
(1- 100(x-0.5)^2)^{1/2}        &   0.4 \leq x \leq 0.6, \\
0                              & \text{otherwise.}
\end{cases}  
\end{eqnarray}
The learned limiters are reported in Fig.~\ref{fig:ensemble_ddl} and analyzed through the accuracy function $\Phi$ (obtained by inverting Eq.~\ref{eq:cphi}). Colors indicate the ranking in terms of shape-preservation performance, from blue (best) to yellow (worst). We find that data-driven limiters that perform worse than the well-know OSTVD3 were excessively diffusive, damping extrema and smoothing sharp features of the solution.

\section{Simplified DDL formulation} \label{sect:appendix-E}

To reduce computational overhead when performing flux-limiting, we perform a model-order reduction by constructing a parametric representation of DDL based on hinge functions, $h(z)=\max(z,0)$, with CFL-dependent knot locations and polynomial amplitudes. The coefficients are obtained with least-squares regression of
\begin{eqnarray}
    y (\delta_c, \nu) &=& \sum_m \sum_k (1-\nu)^k  (\alpha_{m,k} h(\delta_c - c_m(\nu)) + \beta_{m,k} h(c_m(\nu)-\delta_c )) + \sum_k \gamma_k (1-\nu)^k \delta_c \\
c_m(\nu) &=& \sum_j a_{m,j} (1-\nu)^j
\end{eqnarray}
where $\delta_c$ is the curvature, $\nu$ the CFL number and $\alpha_{m,k}$, $\beta_{m,k}$, $\gamma_k$, and $a_{m,j}$ are the regression coefficients (listed in Table~\ref{tab:mor_ddl}). The regression is performed for up-gradient configurations (positive signed curvature), the full expression is recovered by enforcing the symmetry condition described in Section~\ref{sect:mirror}.

\begin{table}[]
    \centering
    \caption{Regression coefficients defining the reduced-order parametric approximation of DDL, including hinge-function amplitudes ($\alpha_{m,k}, \beta_{m,k}$), linear contributions ($\gamma_k$), and CFL-dependent knot locations ($a_{m,j}$).}
    \begin{tabular}{llllll}
Coefficient & Value & Coefficient & Value & Coefficient & Value \\
$a_{0,0}$ & -1.839 & $\alpha_{0,1}$ & -1.819 $\times 10^{-1}$ & $\beta_{0,1}$ & -1.203 $\times 10^{-1}$ \\
$a_{0,1}$ & 1.203 & $\alpha_{0,2}$ & -5.358 $\times 10^{-1}$ & $\beta_{0,2}$ & -4.168 $\times 10^{-1}$ \\
$a_{0,2}$ & -3.746 $\times 10^{-1}$ &  &  &  &  \\
 &  & $\alpha_{1,1}$ & 6.821 $\times 10^{-1}$ & $\beta_{1,1}$ & 7.144 $\times 10^{-1}$ \\
$a_{1,0}$ & -1.060 & $\alpha_{1,2}$ & -1.427 & $\beta_{1,2}$ & -1.547 \\
$a_{1,1}$ & -1.252 $\times 10^{-1}$ &  &  &  &  \\
$a_{1,2}$ & 1.084 $\times 10^{-1}$ & $\alpha_{2,1}$ & 7.828 $\times 10^{-2}$ & $\beta_{2,1}$ & 3.413 $\times 10^{-2}$ \\
 &  & $\alpha_{2,2}$ & 1.445 & $\beta_{2,2}$ & 1.574 \\
$a_{2,0}$ & -9.930 $\times 10^{-1}$ &  &  &  &  \\
$a_{2,1}$ & -8.291 $\times 10^{-2}$ & $\alpha_{3,1}$ & -5.512 $\times 10^{-1}$ & $\beta_{3,1}$ & -6.537 $\times 10^{-1}$ \\
$a_{2,2}$ & 8.796 $\times 10^{-2}$ & $\alpha_{3,2}$ & 4.918 $\times 10^{-1}$ & $\beta_{3,2}$ & 4.217 $\times 10^{-1}$ \\
 &  &  &  &  &  \\
$a_{3,0}$ & -8.258 $\times 10^{-1}$ & $\alpha_{4,1}$ & 2.764 $\times 10^{-1}$ & $\beta_{4,1}$ & 1.996 $\times 10^{-1}$ \\
$a_{3,1}$ & 2.884 $\times 10^{-1}$ & $\alpha_{4,2}$ & -3.432 $\times 10^{-1}$ & $\beta_{4,2}$ & -2.509 $\times 10^{-1}$ \\
$a_{3,2}$ & 7.062 $\times 10^{-2}$ &  &  &  &  \\
 &  & $\alpha_{5,1}$ & -1.982 $\times 10^{-1}$ & $\beta_{5,1}$ & -10.000 $\times 10^{-2}$ \\
$a_{4,0}$ & 1.066 & $\alpha_{5,2}$ & -6.196 $\times 10^{-1}$ & $\beta_{5,2}$ & -6.690 $\times 10^{-1}$ \\
$a_{4,1}$ & -1.329 &  &  &  &  \\
$a_{4,2}$ & 1.011 & $\alpha_{6,1}$ & 1.073 & $\beta_{6,1}$ & 1.123 \\
 &  & $\alpha_{6,2}$ & -1.974 & $\beta_{6,2}$ & -2.028 \\
$a_{5,0}$ & 5.801 $\times 10^{-1}$ &  &  &  &  \\
$a_{5,1}$ & -3.374 $\times 10^{-2}$ & $\alpha_{7,1}$ & -7.184 $\times 10^{-1}$ & $\beta_{7,1}$ & -6.418 $\times 10^{-1}$ \\
$a_{5,2}$ & 4.135 $\times 10^{-1}$ & $\alpha_{7,2}$ & 2.631 & $\beta_{7,2}$ & 2.584 \\
 &  &  &  &  &  \\
$a_{6,0}$ & 9.744 $\times 10^{-1}$ & $\alpha_{8,1}$ & -4.704 $\times 10^{-1}$ & $\beta_{8,1}$ & -5.439 $\times 10^{-1}$ \\
$a_{6,1}$ & 1.132 $\times 10^{-1}$ & $\alpha_{8,2}$ & 3.020 $\times 10^{-1}$ & $\beta_{8,2}$ & 3.731 $\times 10^{-1}$ \\
$a_{6,2}$ & -1.971 $\times 10^{-1}$ &  &  &  &  \\
 &  &  &  &  &  \\
$a_{7,0}$ & 9.771 $\times 10^{-1}$ &  &  &  &  \\
$a_{7,1}$ & 3.010 $\times 10^{-2}$ &  &  &  &  \\
$a_{7,2}$ & 5.463 $\times 10^{-4}$ &  &  &  &  \\
 &  &  &  &  &  \\
$a_{8,0}$ & 1.147 &  &  &  &  \\
$a_{8,1}$ & 5.295 $\times 10^{-2}$ &  &  &  &  \\
$a_{8,2}$ & 7.120 $\times 10^{-1}$ &  &  &  &  \\
 &  &  &  &  &  \\
$\gamma_{1}$ & 1.454 $\times 10^{-2}$ &  &  &  &  \\
$\gamma_{2}$ & 2.565 $\times 10^{-2}$ &  &  &  &  \\
\end{tabular}
    \label{tab:mor_ddl}
\end{table}

\printcredits

\bibliographystyle{model1-num-names} 

\bibliography{bibliography}

\end{document}

%% file: figures/fig_stencil_curv2.tex
\begin{tikzpicture}[scale=1.0]
\newcommand{\offsetplot}[3]{\begin{scope}[shift={(#1,#2)}] #3\end{scope}} 

\newcommand{\xlev}{ 0.}
\newcommand{\alev}{ 0.5}
\newcommand{\blev}{ 1.}
\newcommand{\clev}{ 1.5}

\newcommand{\scale}{0.5} 
\newcommand{\step}{3*\scale}
\newcommand{\half}{0.8*\scale}  
\newcommand{\drawxyz}[4]{%
  \begin{scope}[scale=\scale, every node/.style={scale=\scale, transform shape}]
    \draw[dotted] (-1,#1) -- (0,#2) -- (1,#3);
    \draw (-1-\half,#1) -- (-1+\half,#1);
    \draw (  -\half,#2) -- (  +\half,#2);
    \draw ( 1-\half,#3) -- ( 1+\half,#3);
    
    \ifnum#4=0\relax
      \filldraw[black](0,#2) circle (2pt) (1,#3) circle (2pt);
      \node[rectangle, draw, orange, fill=white] at (-1,#1) {};
    \fi
    \ifnum#4=1\relax
      \filldraw[black] (-1,#1) circle (2pt) (1,#3) circle (2pt);
      \node[rectangle, draw, orange, fill=white] at (0,#2) {};
    \fi
    \ifnum#4=2\relax
      \filldraw[black] (-1,#1) circle (2pt) (0,#2) circle (2pt) ;
      \node[rectangle, draw, orange, fill=white] at (1,#3) {};
    \fi
    \ifnum#4=-1\relax
    \filldraw[black] (-1,#1) circle (2pt) (0,#2) circle (2pt) (1,#3) circle (2pt);
    \fi
  \end{scope}%
}
\newcommand{\dtickq}[2]{
    \draw (#1,\xlev-.1)--(#1,\xlev+.1); \node[below] at (#1,\xlev-.1) {#2}; 
    \draw [dotted] (#1,\xlev-.1)--(#1,\blev-0.2);
    }

\newcommand{\ttickq}[2]{
    \draw (#1,\xlev-.1)--(#1,\xlev+.1); \node[below] at (#1,\xlev-.1) {#2}; 
    }

\newcommand{\drawYAxislabel}[2]{\begin{scope}
    \draw[->, thick] (#1,#2-0.2) -- (#1,#2+0.7) ; \node[left] at (#1-.5,#2+.25) {$\tilde q$} ; 
    \draw (#1-.1,#2) -- (#1+.1,#2); \node[left] at (#1-0.2,#2) {0}; 
    \draw (#1-.1,#2+.5) -- (#1+.1,#2+.5); \node[left] at (#1-0.2,#2+.5) {1}; 
\end{scope}}


\draw[->, thick] (-4*\step-\half,\xlev)--(4*\step+\half,\xlev) node[right]{$\delta^2 \tilde q$};

\offsetplot{-4*\step}{\blev}{ \drawxyz{0. }{1. }{0. }{-1} };
\offsetplot{ 4*\step}{\blev}{ \drawxyz{1. }{0. }{1. }{-1} };

\offsetplot{-3*\step}{\clev}{ \drawxyz{0. }{1. }{0.5}{2} };
\offsetplot{-2*\step}{\clev}{ \drawxyz{0. }{1. }{1. }{-1} };
\offsetplot{-1*\step}{\clev}{ \drawxyz{0. }{0.7}{1. }{1} };
\offsetplot{ 0*\step}{\clev}{ \drawxyz{0. }{0.5}{1. }{1} }; 
\offsetplot{ 1*\step}{\clev}{ \drawxyz{0. }{0.3}{1. }{1} };
\offsetplot{ 2*\step}{\clev}{ \drawxyz{0. }{0. }{1. }{-1} };
\offsetplot{ 3*\step}{\clev}{ \drawxyz{0.5}{0. }{1.}{0} };

\offsetplot{-3*\step}{\alev}{ \drawxyz{0.5}{1. }{0. }{0} };
\offsetplot{-2*\step}{\alev}{ \drawxyz{1. }{1. }{0. }{-1} };
\offsetplot{-1*\step}{\alev}{ \drawxyz{1. }{0.7}{0. }{1} };
\offsetplot{ 0*\step}{\alev}{ \drawxyz{1. }{0.5}{0. }{1} }; 
\offsetplot{ 1*\step}{\alev}{ \drawxyz{1. }{0.3}{0. }{1} };
\offsetplot{ 2*\step}{\alev}{ \drawxyz{1. }{0. }{0. }{-1} };
\offsetplot{ 3*\step}{\alev}{ \drawxyz{1. }{0. }{0.5}{2} };

\ttickq{-4*\step}{$-2$}
\ttickq{-2*\step}{$-1$}
\ttickq{ 0*\step}{$0$}
\ttickq{ 2*\step}{$1$}
\ttickq{ 4*\step}{$2$}

\drawYAxislabel{-4.6*\step}{\blev}


\end{tikzpicture}